\documentclass[10pt]{amsart}
\usepackage{amsmath,amssymb,epsf, epic, eepic}
\usepackage{psfig}
\usepackage{amsfonts}
\usepackage{times}

\input xy
\xyoption{all}

\newtheorem{thm}{Theorem}[section]
\newtheorem{lemma}[thm]{Lemma}
\newtheorem{corollary}[thm]{Corollary}
\newtheorem{prop}[thm]{Proposition}
\newtheorem{exe}[thm]{Example}

\theoremstyle{definition}

\newcommand{\bN}{{\mathbb{N}}}
\newcommand{\bZ}{{\mathbb{Z}}}

\newcommand{\bF}{{\mathbb{F}}}
\newcommand{\bK}{{\mathbb{K}}}

\newcommand{\cA}{{\mathcal{A}}}
\newcommand{\cB}{{\mathcal{B}}}
\newcommand{\cC}{{\mathcal{C}}}

\newcommand{\cK}{{\mathcal{K}}}

\DeclareMathOperator{\Ker}{Ker}
\DeclareMathOperator{\Image}{Im}

\newcommand{\ra}{\rightarrow}
\newcommand{\op}[1]{\overline{#1}}

\newcommand{\spaces}{\;\;\;\;\;\;\;}
\newcommand{\ot}{\otimes}
\newcommand{\ol}{\overline}
\newcommand{\ob}{\overline}

\newcommand{\FF}{\bF_{2}}
\newcommand{\FFu}{\bF_2[u]}

\newcommand{\kh}[3]{\mbox{KH}^{{#1},{#2}}_{\FF}({#3})}
\newcommand{\rkh}[3]{\widetilde{\mbox{KH}}^{{#1},{#2}}_{\FF}({#3})}

\newcommand{\bn}[3]{\mbox{BN}^{{#1},{#2}}({#3})}

\newcommand{\sbn}[2]{\mbox{BN}_s^{#1}({#2})}

\newcommand{\kk}[3]{\bK^{{#1},{#2}}({#3})}

\newcommand{\khch}[3]{\ob{\cC}^{{#1},{#2}}({#3})}

\newcommand{\ch}[2]{C^{{#1},{#2}}}
\newcommand{\cchb}[2]{\op{\cC}^{{#1},{#2}}}
\newcommand{\ccht}[2]{\widetilde{\cC}^{{#1},{#2}}}
\newcommand{\cch}[2]{\cC^{{#1},{#2}}}
\newcommand{\scchb}[1]{\op{\cC}^{{#1}}}
\newcommand{\sccht}[1]{\widetilde{\cC}^{{#1}}}
\newcommand{\scch}[1]{\cC^{{#1}}}
\newcommand{\ubn}[1]{\mbox{BN}^{{#1}}(L)^\prime}
\newcommand{\rubn}[1]{\widetilde{\mbox{BN}}^{{#1}}(L)^\prime}

\newcommand{\dd}[2]{\partial_{#1}^{#2}}
\newcommand{\bb}[2]{\beta_{#1}^{#2}}

\newcommand{\cd}[1]{\cC (D(*{#1}))}

\newcommand{\cdt}[1]{\cC (\tilde D(*{#1}))}

\newcommand{\lk}{\mbox{{\em lk}}}

\title{Calculating Bar-Natan's characteristic
  two  Khovanov homology}
\author{Paul Turner}
\address{School of Mathematical and Computer Sciences \\Heriot-Watt
  University\\ Edinburgh EH14 4AS\\Scotland}
\email{paul@ma.hw.ac.uk}
\begin{document}


\begin{abstract}
We investigate Bar-Natan's characteristic
two Khovanov link homology theory studying both the filtered and
bi-graded theories. The filtered theory is computed explicitly and the
bi-graded theory analysed by setting up a family of spectral sequences. The 
$E_2$-pages can be described in terms of groups arising from the
action of a certain endomorphism on $\FF$-Khovanov homology. Some
simple consequences are discussed.
\end{abstract}

\maketitle

\section{Introduction and statement of results}

In his remarkable paper \cite{barnatan2} Bar-Natan shows that any
Frobenius algebra satisfying certain conditions gives rise to a
homology theory in the sense of Khovanov's homology for links
\cite{khovanov}. In particular he singles out a characteristic two
theory which associates to each link diagram a bi-graded
$\FF[u]$-module in such a way that any two diagrams representing the
same link give isomorphic modules. By setting $u=1$ one loses the
bi-grading obtaining a singly graded theory with a filtration in place
of the internal grading. The purpose of the current paper is to
investigate how to calculate both the filtered and bi-graded versions
of Bar-Natan's characteristic two link homology theory.

In order to establish some notation let us briefly recall some things
about Khovanov's homology. Let $L$ be an oriented link and $D$ a
diagram for $L$. Let $\khch * * D$ be Khovanov's complex over $\FF$
for the diagram $D$ and let $\kh * * L$ be the resulting
$\FF$-Khovanov homology \cite{khovanov} (see also
\cite{barnatan1}). The differential in $\khch **D$ will be denoted
$\partial\colon \khch i j D \ra \khch {i+1} j D$. For $v$ belonging to
the vector space associated to the complete smoothing $\alpha$ this is
defined by
\[
\partial (v) = \sum  \cA (S_{e})(v)
\]
where the summation is over all edges $e$ in the cube of smoothings
whose tail is $\alpha$ and $S_e$ is the cobordism attached to that
edge. The signs drop out because we are working mod 2. $\cA$ is the
TQFT associated to the Frobenius algebra $A=\FF \{1,x\}$ with
multiplication $\ob m$ given by
\[
\ob m(1,1) = 1 \spaces \ob m(1,x) = x \spaces \ob m(x,1)=x \spaces \ob m(x,x) = 0
\]
and comultiplication $\ob \Delta$ given by
\[
\ob \Delta(1) = 1\ot x + x\ot 1  \spaces \ob \Delta(x) = x\ot x.
\] 
The unit and counit are given by
\[
i(1) = 1 \spaces \epsilon (1) =0 \spaces \epsilon (x) =1.
\] 

Khovanov homology is bi-graded and we refer to the first grading as the {\em homological grading} and the second grading as the $q$-grading. Given $v\in\kh **L$ we write $q(v)$ for its $q$-grading.

Bar-Natan's theory is defined over $\FFu$
where the bi-degree of $u$ is $(0,-2)$. This theory is defined analogously to
Khovanov's using the TQFT associated to the Frobenius algebra $\FFu
\{1,x\}$ with multiplication $ m$ given by
\[
m (1,1) = 1 \spaces m (1, x) = x \spaces m (x, 1)=x
\spaces m (x, x) = ux
\]
comultiplication $\Delta$
\[
\Delta(1) = 1\ot x + x\ot 1 + u 1\ot 1 \spaces \Delta(x) = x\ot x
\] 
and unit and counit
\[
i(1) = 1 \spaces \epsilon (1) =0 \spaces \epsilon (x) =1.
\] 
The complex associated to the diagram $D$ will be denoted $\cch ** (D)$ with differential $d$. We denote the resulting bi-graded homology of the link $L$ by $\bn ** L$. This will be referred to as {\em graded Bar-Natan theory}.

By setting $u=1$ the underlying groups are those of Khovanov theory but the differential $d$ no longer respects the $q$-grading. Thus ignoring the second grading we can consider $\scchb * (D)$ as a complex with differential $d$. (We will suppress from the notation the fact that $d$ now has $u$ set to $1$). In fact while the second grading is not preserved under $d$ it does not decrease which gives rise to a filtration on $(\scchb * (D),d)$. We denote the resulting link homology by $\ubn * $ referring to this as {\em filtered Bar-Natan theory}.

In this paper we define, in Section~2, an endomorphism $\beta_*\colon
\kh * * L \ra \kh {*} {*} L$ of bi-degree $(1,2)$ on $\FF$-Khovanov
homology. Since $\beta_*^2=0$ this means $\beta_*$ can be viewed as a
differential on $\kh ** L$ and we can take homology to obtain
secondary groups $\kk * * L$. These secondary groups appear later in
the $E_2$-page of certain spectral sequences. In Section~3 we study
filtered Bar-Natan theory of an oriented link $L$. Following the work of Lee \cite{lee} we
calculate this explicitly and the first result is as follows.

\vspace*{6pt}

{\bf Theorem }\ref{thm:dim}
{\em
The dimension of $\ubn * $ is $2^k$ where $k$ is the number of
components in $L$. Moreover if $L_1, \ldots , L_k$ are the components
then
\[
\dim (\ubn i )= \mbox{Card} \{ E\subset \{1,2, \cdots , k \} \mid
2 \sum_{l\in E,m\notin E} \lk(L_l,L_m) = i\}
\]
where  $\lk(L_l,L_m)$ is the linking number between component $L_l$
and $L_m$.
} 

\vspace*{6pt}

The filtration of the filtered theory gives rise to a spectral sequence whose $E_2$-page can be identified.

\vspace*{6pt}

{\bf Theorem }\ref{thm:ss1}
{\em There is a spectral sequence with $E_1$-page $\kh ** L$ converging to
$\ubn *$. The $E_2$-page is given in terms of the secondary groups $\kk ** L$.}

\vspace*{6pt}

We end this section by showing that these results carry over to the setting of reduced link homology theory.

In Section~4 we study the graded theory, where there is a spectral sequence for each $q$-grading.

\vspace*{6pt}

{\bf Theorem }\ref{thm:ss}
{\em Given $j\in\bZ$ there is a spectral sequence with $E_1$-page determined by the mod 2 Khovanov homology of $L$ converging to $\bn * j L$. The $E_2$-page can be determined using the endomorphism $\beta_*$.}

\vspace*{6pt}

We then show that graded Bar-Natan theory stabilises with respect to the $q$-grading and prove that the singly graded stable limit is isomorphic to the singly graded (filtered) theory. In fact the spectral sequence of Theorem \ref{thm:ss} coincides with the spectral sequence of Theorem \ref{thm:ss1} for $j$ in the stable range.

 In Section~5 we discuss two elementary consequences of
these results.  We give a brief
discussion of the appearance of ``shifted pawn moves'' in graded Bar-Natan
theory and prove the following result about the form of the mod 2 Khovanov polynomial for thin knots.  

\vspace*{6pt}

{\bf Theorem }\ref{thm:thin}
{\em If $L$ is an $\FF$H-thin knot with homology concentrated on
diagonals $j=s-1+2i$ and $j=s+1+2i$ then there exists a polynomial
$Kh^\prime (L)$ such that 
\[
Kh_{\FF}(L) = q^{s-1}(1+q^2)(1+(1+tq^2)Kh^\prime (L))
\]
where  $Kh^\prime (L)$ is a polynomial in $tq^2$.}

\vspace*{6pt}

Throughout we draw heavily on the techniques developed by Lee in
\cite{lee} and in Section~6 we include for completeness the details of
the proof that the endomorphism $\beta_*$ is independent of the chosen
diagram.


\label{auto}\section{An endomorphism on $\FF$-Khovanov homology}

We will define an endomorphism 
$
\beta_*\colon \kh * * L \ra \kh {*}
{*} L
$
 of bi-degree $(1,2)$ in a similar manner to the map $\Phi$ in rational Khovanov homology defined by
Lee in \cite{lee}. Define operations $\tilde m \colon A \otimes A \ra A$ and $\tilde
\Delta \colon A \ra A \otimes A $  by
\[
\tilde m (1, 1) =0 \spaces \tilde m (1,x)=0 \spaces \tilde m ( x, 1)  = 0
\spaces \tilde m ( x, x) =x
\]
and
\[
\tilde \Delta(1) = 1\ot 1 \spaces \tilde \Delta (x) = 0.
\]
There is no compatible unit or counit so there is no Frobenius algebra
structure. None the less given a diagram $D$ and an edge $e$ of the cube of $D$ we can
define a map $\cB(S_e)$ for the
cobordism $S_e$ associated to  $e$. This is possible since such a
cobordism is made up of cylinders and a pair-of-pants surface only. This 
allows us to define a map 
$
\beta\colon \khch i j D \ra \khch {i+1}
{j+2} D
$
 in an analogous fashion to the differential: for $v$
belonging to the vector space associated to the complete smoothing
$\alpha$ set 
\[
\beta(v) = \sum \cB(S_e)(v).
\]
where as before the sum is over all edges whose tail is $\alpha$.

\begin{lemma} {\it The map $\beta$ is a map of complexes and  $\beta^2=0$.}
\end{lemma}  

\begin{proof}
The proof is identical to Lee's proof in \cite{lee} for her map
$\Phi$. For the first part it suffices to show the following three
equations.
\begin{eqnarray*}
\ol m\circ (\tilde m \ot 1) + \tilde m \circ (\ol m \ot 1) + \ol m \circ (1\ot \tilde m) +
\tilde m \circ (1\ot \ol m) =0 & &\\
(\ol \Delta \otimes 1)\circ \tilde \Delta + (\tilde \Delta \otimes 1)\circ
\ol \Delta + (1\ot \ol \Delta)\ot \tilde \Delta +  (1\ot \tilde \Delta)\ot \ol \Delta =0
&&
\\
\ol \Delta \circ \tilde m + \tilde \Delta \circ \ol m + (\ol m\ot 1)\circ (1\ot \tilde
\Delta) + (\tilde m\ot 1)\circ (1\ot \ol \Delta) = 0 && 
\end{eqnarray*}
These can almost instantly be verified.

The second part follows straightforwardly from the fact that $\tilde m$ (resp. $\tilde
\Delta$)  is
(co)commutative and (co)associative and that
$
\tilde \Delta \circ \tilde m = (\tilde m \ot 1)\circ (1\ot \tilde \Delta)
$ 
all of which are again simply verified.
\end{proof}

This means that $\beta$ induces a map $\beta_*$ in homology, however
{\em a priori} $\beta_*$ depends on the diagram $D$.  It turns out
there is no such dependence and any diagram for $L$ gives the same
induced map in homology. Or, more precisely, $\beta_*$ commutes with
the isomorphisms in homology induced by Khovanov's quasi-isomorphisms
of complexes given by Reidemeister moves.  This follows from the
following proposition. (We chose to write Reidemeister I with
negative twist as a combination of other moves.)


\begin{prop}\label{prop:commute}
{\it The map $\beta$ commutes up to boundaries with Khovanov's quasi-isomorphisms
for Reidemeister I positive twist, Reidemeister II and Reidemeister III.}
\end{prop}

Again the proof follows the work of Lee \cite{lee} almost verbatim. For
completeness we give the details in Section 6. 

Thus, given an oriented link $L$ there is a well defined map
\[
\beta_*\colon \kh i j L \ra \kh {i+1} {j+2} L.
\]
The action of $\beta_*$ on the $\FF$-Khovanov homology is additional
information about the link.

\begin{exe}\label{exe:trefoil}
{\em  In this example  we compute the action of $\beta_*$ when $L=$ {\raisebox{-3mm}{\psfig{figure=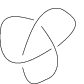}}}. The $\FF$-Khovanov homology is
  summarized in Figure \ref{fig:khtref}.

\begin{figure}[h]
\begin{center}
\setlength{\unitlength}{0.00033333in}
\begingroup\makeatletter\ifx\SetFigFont\undefined%
\gdef\SetFigFont#1#2#3#4#5{%
  \reset@font\fontsize{#1}{#2pt}%
  \fontfamily{#3}\fontseries{#4}\fontshape{#5}%
  \selectfont}%
\fi\endgroup%
{\renewcommand{\dashlinestretch}{30}
\begin{picture}(3387,4062)(0,-10)
\path(975,3012)(3375,3012)(3375,12)
	(975,12)(975,3012)
\path(975,2412)(3375,2412)
\path(975,1812)(3375,1812)
\path(3375,1212)(975,1212)
\path(975,612)(3375,612)
\path(1575,3012)(1575,12)
\path(2175,3012)(2175,12)
\path(2775,3012)(2775,12)
\put(2925,3162){\makebox(0,0)[lb]{{\SetFigFont{8}{9.6}{\rmdefault}{\mddefault}{\updefault}0}}}
\put(2400,3162){\makebox(0,0)[lb]{{\SetFigFont{8}{9.6}{\rmdefault}{\mddefault}{\updefault}-1}}}
\put(1200,3162){\makebox(0,0)[lb]{{\SetFigFont{8}{9.6}{\rmdefault}{\mddefault}{\updefault}-3}}}
\put(600,2562){\makebox(0,0)[lb]{{\SetFigFont{8}{9.6}{\rmdefault}{\mddefault}{\updefault}-1}}}
\put(600,1887){\makebox(0,0)[lb]{{\SetFigFont{8}{9.6}{\rmdefault}{\mddefault}{\updefault}-3}}}
\put(600,1362){\makebox(0,0)[lb]{{\SetFigFont{8}{9.6}{\rmdefault}{\mddefault}{\updefault}-5}}}
\put(600,762){\makebox(0,0)[lb]{{\SetFigFont{8}{9.6}{\rmdefault}{\mddefault}{\updefault}-7}}}
\put(1725,3162){\makebox(0,0)[lb]{{\SetFigFont{8}{9.6}{\rmdefault}{\mddefault}{\updefault}-2}}}
\put(1125,237){\makebox(0,0)[lb]{{\SetFigFont{8}{9.6}{\rmdefault}{\mddefault}{\updefault}1}}}
\put(1125,837){\makebox(0,0)[lb]{{\SetFigFont{8}{9.6}{\rmdefault}{\mddefault}{\updefault}1}}}
\put(1800,837){\makebox(0,0)[lb]{{\SetFigFont{8}{9.6}{\rmdefault}{\mddefault}{\updefault}1}}}
\put(1800,1437){\makebox(0,0)[lb]{{\SetFigFont{8}{9.6}{\rmdefault}{\mddefault}{\updefault}1}}}
\put(2925,2037){\makebox(0,0)[lb]{{\SetFigFont{8}{9.6}{\rmdefault}{\mddefault}{\updefault}1}}}
\put(2925,2637){\makebox(0,0)[lb]{{\SetFigFont{8}{9.6}{\rmdefault}{\mddefault}{\updefault}1}}}
\put(0,1437){\makebox(0,0)[lb]{{\SetFigFont{8}{9.6}{\rmdefault}{\mddefault}{\updefault}j}}}
\put(1950,3837){\makebox(0,0)[lb]{{\SetFigFont{8}{9.6}{\rmdefault}{\mddefault}{\updefault}i}}}
\put(600,162){\makebox(0,0)[lb]{{\SetFigFont{8}{9.6}{\rmdefault}{\mddefault}{\updefault}-9}}}
\end{picture}
}
\end{center}
\caption{$\FF$-Khovanov homology of the trefoil}
\label{fig:khtref}
\end{figure}

For dimensional reasons there are only two potentially non-zero maps
to consider, namely
\[
\beta_*\colon \kh {-3} {-9} L \ra \kh {-2} {-7}L
\]
and
\[
\beta_*\colon \kh {-3} {-7} L \ra \kh {-2} {-5}L.
\]
A generating cycle for $\kh {-3} {-9} L $ is given by $x\otimes x
\otimes x$ and noting that $\khch {-2} {*} L = A^{\ot 2}
\oplus A^{\ot 2} \oplus A^{\ot 2}  $ we calculate 
\[
\beta (x\otimes x \otimes x) = x\otimes x
+  x^\prime\otimes x^\prime +  x^{\prime\prime}\otimes x^{\prime\prime} \in \khch {-2} {*} L
\]
where the primes indicate different copies of $A^{\ot 2}$.
Furthermore, by looking at the
$\FF$-Khovanov differential $\partial$  the element on the
right is easily seen to be a generator and so $ \beta_*\colon \kh {-3}
{-9} L \ra \kh {-2} {-7}L$ is an isomorphism.

Similarly  $1\ot x \ot x + x \ot 1 \ot x + x\ot x
\ot 1 $ is a generator of $\kh {-3} {-7} L $ and
\[
\beta(1\ot x \ot x + x \ot 1 \ot x + x\ot x
\ot 1) = x\ot 1 +  1\ot x^\prime +  x^{\prime\prime} \otimes 1\in \khch {-2} {*} L 
\]
which is  easily seen to be a generator of $\kh {-2} {-7} L$.  Thus
$\beta_*\colon \kh {-3} {-7} L \ra \kh {-2} {-5}L $ is an isomorphism
too.}
\end{exe}

The action of $\beta_*$ may be encoded by
regarding it as a differential on $\FF$-Khovanov homology and
considering the associated homology groups. More precisely for each
$k$ define a complex $\cK(k)^*$ by $\cK(k)^i = \kh i {k+2i} L
$ with differential $\beta_*$. Set
\[
\kk i j L = H^i(\cK(j-2i)^*, \beta_*)
= 
\frac{\Ker (\beta_*\colon \kh i j L \ra \kh {i+1} {j+2} L)}
{\Image (\beta_*\colon \kh {i-1} {j-2} L \ra \kh {i} {j} L)}.
\]
From the above we know that these groups are link invariants and they
contain secondary information about the link. One can take the
associated Poincar\'e polynomial to get a {\em secondary link polynomial},
namely
\[
P(t,q)(L)= \sum_{i,j\in \bZ} t^i q^j \dim_{\FF}(\kk i j L ).
\]

\begin{exe}
{\em Using the computation above for the trefoil we see 
\[
P(t,q)( {\raisebox{-3mm}{\psfig{figure=trefoil.eps}}}) = \frac{1}{q} + \frac 1 {q^3}.
\]}
\end{exe}


\section{Filtered Bar-Natan theory: setting $u=1$}\label{sec:filtered}

By setting $u=1$ in Bar-Natan theory the vector spaces in $\scch *$
transform into the vector spaces in $\scchb *$. While Bar-Natan's
differential $d$ (setting $u=1$) does not respect the second grading
in $\scchb *$ is it easy to see that it cannot decrease this grading
and so a filtration on $\scchb *$ can be defined. This is analogous
the to situation in Lee's theory \cite{lee} (see also
\cite{rasmussen}) and with minor adjustments the work of Lee
carries over to this situation. We denote the singly graded $u=1$
Bar-Natan theory by $\ubn * $. In this section we show that Lee's
methods to compute her theory explicitly can be applied to $\ubn *
$. We also analyse the spectral sequence arising from the filtration
and show that in our case, in a departure from Lee's work, we can work
with {\em reduced} homology as well.

\subsection{Explicit calculation of $\ubn *$}
Lee has a clever argument for calculating her theory explicitly,
which we adapt to compute $\ubn *$.

\vspace*{6pt}

\begin{thm}\label{thm:dim} 
{\it The dimension of $\ubn * $ is $2^k$ where $k$ is the number of
components in $L$. Moreover if $L_1, \ldots , L_k$ are the components
then
\[
\dim (\ubn i ) = \mbox{Card} \{ E\subset \{1,2, \cdots , k \} \mid
2 \sum_{l\in E,m\notin E} \lk(L_l,L_m) = i\}
\]
where  $\lk(L_l,L_m)$ is the linking number between component $L_l$
and $L_m$.}
\end{thm}

\vspace*{6pt}

To prove this theorem we first perform a change of basis on the
Frobenius algebra $A$ by setting
\[
a = x+1 \;\;\;\; \text{ and } \;\;\;\; b = x.
\]
The new basis $\{a,b\}$ diagonalises the Frobenius structure and we have:
\[
aa = a \spaces ab=ba=0 \spaces bb=b
\]
with comultiplication 
\[
a \mapsto a\ot a \spaces b\mapsto b\ot b 
\] 
and unit and counit
\[
i(1) = a+b \spaces \epsilon (a) = \epsilon (b)  = 1.
\] 
Given an orientation of $L$ one can choose the unique orientation-preserving complete smoothing and partition the circles of such a smoothing
into two groups, Group A and Group B (see \cite{lee} or
\cite{rasmussen}). By assigning the element $a$ to those circles in Group A
and $b$ to those circles in Group B we obtain an element of $\scchb * (L)$
and it is easy to see this element is a cycle for the differential $d$. As there are $2^k$ possible orientations of $L$ this gives $2^k$ cycles.

Using the fact that $\scchb *$ is constructed from a Frobenius algebra
there is an inner product on $\scchb *$ from which we can define an
adjoint differential $d^*$ and a simple computation shows that the
$2^k$ cycles above are also cycles with respect to this adjoint
differential. Using the duality properties of this adjoint one can
identify $\ubn *$ with $ \{ z\in \scchb i (D) \mid d z = 0 \mbox{ and } d^*
z = 0 \}$. From here on Lee's argument given in \cite{lee} works almost
verbatim showing that the $2^k$ classes above generate the homology
and moreover that their homological degree is as given in the theorem.

\subsection{A spectral sequence}\label{sec:filteredss}
An arbitrary element of $\scchb i$ is not homogeneous with respect to
the second grading but may be written as a sum $\sum_{\lambda\in
\Lambda}v_\lambda $ of homogeneous elements $v_\lambda$
 for some indexing set $\Lambda$. Let $\gamma$ be the mod 2 number of components of the link and recall that $\cchb i j =0$ unless $j=\gamma$ mod 2. We filter
$\scchb *$ by setting
\[
F^k\scchb i = \{\sum_{\lambda\in \Lambda} v_\lambda \in \scchb i \mid q(v_\lambda) \geq 2k+\gamma \text{ for } \lambda\in\Lambda\}.
\]
We have $dF^k\subset F^k$ and since only finitely many
$\cchb i j$ are non-trivial it follows that the filtration is bounded
and so there is a spectral sequence converging to $\ubn
*$. The $E_0$-page of this spectral sequence is given by
\[
E_0^{k,l} = \frac{F^k\scchb {k+l}}{F^{k+1}\scchb {k+l}} = \cchb {k+l}
{2k+\gamma}.
\]
Notice that Bar-Natan's differential $d$ can be written as the sum of Khovanov's differential $\partial$ and the map $\beta$. From this we see  that the zero'th differential in the spectral
sequence $d_0\colon E_0^{k,l} \ra E_0^{k,l+1}$ agrees with $\partial$. Thus the $E_1$-page is given by
\[
E_1^{k,l} = H^{k+l}(\cchb * {2k+\gamma}) = \kh {k+l} {2k+\gamma} L.
\]
The first differential 
$d_1\colon E_1^{k,l} \ra
E_1^{k+1,l}$ is given by the boundary map in the long exact sequence
associated to the following short exact sequence of complexes.
\[
0 \ra \frac{F^{k+1}}{F^{k+2}} \ra  \frac{F^{k}}{F^{k+2}} \ra
\frac{F^{k}}{F^{k+1}} \ra 0 
\]
Since $d=\partial +\beta$ it follows that the boundary map is given by
$x\mapsto \beta_*(x)$ showing that $d_1=\beta_*$. Thus $E_2^{k,l} =
\kk {k+l} {2k+\gamma} L$.

In summary we have:

\vspace*{6pt}

\begin{thm}\label{thm:ss1}
{\it There is a spectral sequence with $E_1$-page $\kh ** L$ converging to
$\ubn *$. The $E_2$-page is given in terms of the secondary groups $\kk ** L$.}
\end{thm}

\subsection{Reduced homology for knots}
The results above carry over to reduced homology and I am grateful to
J. Rasmussen and to the referee for pointing this out to me.  We first
recall the definition of reduced mod 2 Khovanov homology of a knot $L$
(see \cite{khovanov2}). Let $D$ be a diagram of $L$ and choose a
preferred point on $D$ referring to this as the basepoint. Letting $U$
denote the unknot we consider the cobordism from the disjoint union
$U\sqcup D$ to $D$ given by the Morse 1-handle move fusing the unknot
to $D$ at the base point. This induces a map of complexes $A \otimes
\cchb * * (D) \ra \cchb ** (D)$ turning $\cchb ** (D)$ into complex of
$A$-modules. Khovanov argues that the homotopy type of this complex is
an invariant of $L$. Now let $Q=A/xA\cong \FF$ be the one dimensional
representation of $A$ and form the {\em reduced} chain complex by
setting
\[
\ccht ** (D) = \cchb **(D) \otimes_A Q.   
\]
The  {\em reduced mod 2 Khovanov
homology} $\rkh ** L$ is the homology of this complex. Up to isomorphism this does
not depend on the diagram chosen nor on the choice of base point. 

One can easily check that Bar-Natan's differential $d$ (setting $u=1$)
is also well defined on $\sccht * (D)$ though again this does not
respect the second grading. This leads in a similar manner to {\em
reduced filtered Bar-Natan theory} $\rubn * $. Furthermore, the map
$\beta$ of Section~2 induces a map $\tilde \beta$ on
$\ccht ** (D)$. This is very different to the situation with Lee's
theory where her map $\Phi$ does not descend to the reduced rational
theory. The spectral sequence of the previous section can now be
constructed  in the reduced setting too. That is to say by using the
filtration of $(\sccht * (D),d)$ which is induced by the second grading in $\ccht ** (D)$ there is a spectral sequence with
$E_1$-page given by
\[
E_1^{k,l} = H^{k+l}(\ccht * {2k+\gamma+1}) = \rkh {k+l} {2k+\gamma+1} L.
\]
converging to reduced Bar-Natan theory $\rubn *$. The differential
$d_1$ can be identified with $\tilde\beta$.

\section{Graded Bar-Natan theory}
\subsection{Spectral sequences}\label{sec:gradedss}
In this subsection we set up spectral sequences for computing the
graded Bar-Natan theory $\bn **L$ of a link $L$. A different spectral
sequence is needed to compute each fixed $q$-grading of the Bar-Natan
theory, that is, given $j\in\bZ$ there is a spectral sequence which
computes $\bn * j L$. These spectral sequences will be along the lines
of that constructed in Khovanov \cite{khovanov} starting with Khovanov
theory over $\bZ$ converging to Khovanov theory over $\bZ[c]$.

Now we fix $j\in \bZ$ for the rest of the section. To abbreviate
notation a little we will omit $L$ writing $\cch **$ for $\cch **
(L)$. Since $\cch * *$ is defined over $\FF[u]$ we can consider the
multiplication by $u$ map which we will write $u\colon \cch ** \ra
\cch * {*-2}$. We can use this to filter $\cch *j$ by setting
\[
F^k\cch*j = u^k\cch** \cap \cch *j.
\]
Equivalently there is an isomorphism of groups
\[
\cch * j \cong \bigoplus_{p\geq 0} \cchb * {j+2p}
\]
and the filtration on  $\ch * j$ is given by
\[
F^k \cch * j = \bigoplus_{p\geq k} \cchb * {j+2p}.
\]
Since only finitely many groups in $\cchb **$ are non-trivial it
follows that this gives a bounded filtration. Moreover it is easy to
see that $dF^k \subseteq F^k$. Associated to this filtration we get a
spectral sequence which converges to $H^*(\cch * j) = \bn {*} j
L$. There are no problems with convergence since the filtration is
bounded. More precisely there is a filtration on $H^*(\cch * j)$ induced by the
above filtration by setting
\[
F^k\bn{*} j L = \Image (H^*(F^k\cch * j) \ra H^*(\cch * j))
\]
and the spectral sequence has $E_\infty$-term given by
\[
E_\infty^{k,l} = \frac{F^kH^{k+l}(\cch * {j})}{F^{k+1}H^{k+l}(\cch *
  {j})} = 
\frac{F^k\bn {k+l}{j} L}{F^{k+1}\bn {k+l}{j} L}
\]
and so 
\[
\bn {i} j L \cong \bigoplus_{k+l=i}  E_\infty^{k,l}.
\]

We will now identify the $E_1$-page, the differential $d_1$ and hence the
$E_2$-page of this spectral sequence. We will require the following
Lemma.

\begin{lemma} \label{lem:db}
{\it If $v\in \cchb i {j+2p} \subset \cch i j $ then 
\[
d(v) = \partial (v) + \beta_*(v)
 \in \cchb {i+1} {j+2p}\oplus \cchb {i+1} {j+2(p+1)} \subset \cch {i+1} j. 
\]
}
\end{lemma}

\begin{proof}
This is immediate after observing that $m = \ol m + u \tilde m$ and
$\Delta = \ol \Delta + u \tilde \Delta$.
\end{proof}

The $E_0$-page of the spectral sequence is given by
\[
E_0^{k,l} = 
\begin{cases}
\frac{F^k\cch {k+l} j }{F^{k+1}\cch {k+l} j } = \cchb {k+l}
{j+2k} & k\geq 0\\
0 & k<0
\end{cases}
\]
It is immediate from the lemma above that the differential $d_0\colon
E_0^{k,l} \ra E_0^{k,l+1}$ which is induced from $d$ agrees with
$\partial$.  Thus the $E_1$-term is given by
\[
E_1^{k,l} =
\begin{cases}  H^{k+l}(\cchb * {j+2k}) = \kh {k+l} {j+2k} L & k\geq 0\\
0 & k<0
\end{cases}
\]

We now claim that the differential $d_1\colon E_1^{k,l} \ra
E_1^{k+1,l}$ agrees with $\beta_*$ defined in Section~2.
The differential $d_1$ is given by the boundary map in the long exact sequence
associated to the following short exact sequence of complexes:
\[
0 \ra \cchb * {j+2(k+1)} \ra \cchb * {j+2(k+1)} \oplus \cchb * {j+2k} \ra
\cchb * {j+2k} \ra 0.
\]  
Using the lemma above it is easy to see that the boundary map of the
associated long exact sequence is given by $x\mapsto \beta_*(x)$
showing that $d_1=\beta_*$. Thus we can identify the $E_2$-page as follows.
\[
E_2^{k,l} = 
\begin{cases}
\kk {k+l}{j+2k} L & k> 0\\
\Ker (\beta_*\colon \kh {l} j L \ra \kh {l+1} {j+2} L ) & k = 0\\
0 & k<0
\end{cases}
\]

To summarize:

\vspace*{6pt}

\begin{thm}\label{thm:ss}
{\it Given $j\in\bZ$ there is a spectral sequence with $E_1$-page determined by the mod 2 Khovanov homology of $L$ converging to $\bn * j L$. The $E_2$-page can be determined using the endomorphism $\beta_*$.}
\end{thm}

\vspace*{6pt}

\begin{exe}\label{exe:tref}
{\em We now use the spectral sequence to compute the Bar-Natan theory for
the trefoil in Example \ref{exe:trefoil}.  The 
$E_1$-pages for the spectral sequences
corresponding to $j=-1$,$-3$, $-5$, $-7$, $-9$ are given in Figure
\ref{fig:trefE1}.
\begin{figure}[h]
\begin{center}
\setlength{\unitlength}{0.00033333in}
\begingroup\makeatletter\ifx\SetFigFont\undefined%
\gdef\SetFigFont#1#2#3#4#5{%
  \reset@font\fontsize{#1}{#2pt}%
  \fontfamily{#3}\fontseries{#4}\fontshape{#5}%
  \selectfont}%
\fi\endgroup%
{\renewcommand{\dashlinestretch}{30}
\begin{picture}(13053,2727)(0,-10)
\put(30,360){\makebox(0,0)[lb]{{\SetFigFont{5}{6.0}{\rmdefault}{\mddefault}{\updefault}-5}}}
\path(375,1500)(1875,1500)
\path(375,1200)(1875,1200)
\path(375,900)(1875,900)
\path(375,600)(1875,600)
\path(675,1800)(675,300)
\path(975,1800)(975,300)
\path(1275,1800)(1275,300)
\path(1575,1800)(1575,300)
\path(375,300)(375,900)
\path(375,900)(375,2700)
\path(405.000,2580.000)(375.000,2700.000)(345.000,2580.000)
\path(375,2100)(1875,2100)(1875,1800)
\path(375,2400)(1875,2400)(1875,2100)
\path(975,2400)(975,1800)
\path(675,2400)(675,1800)
\path(1275,2400)(1275,1800)
\path(1575,2400)(1575,1800)
\put(2100,0){\makebox(0,0)[lb]{{\SetFigFont{8}{9.6}{\rmdefault}{\mddefault}{\updefault}$k$}}}
\path(375,300)(2175,300)
\path(2055.000,270.000)(2175.000,300.000)(2055.000,330.000)
\put(150,2400){\makebox(0,0)[lb]{{\SetFigFont{8}{9.6}{\rmdefault}{\mddefault}{\updefault}$l$}}}
\path(3075,1800)(4575,1800)(4575,300)
\path(3075,1500)(4575,1500)
\path(3075,1200)(4575,1200)
\path(3075,900)(4575,900)
\path(3075,600)(4575,600)
\path(3375,1800)(3375,300)
\path(3675,1800)(3675,300)
\path(3975,1800)(3975,300)
\path(4275,1800)(4275,300)
\path(3075,300)(3075,900)
\path(3075,900)(3075,2700)
\path(3105.000,2580.000)(3075.000,2700.000)(3045.000,2580.000)
\path(3075,2100)(4575,2100)(4575,1800)
\path(3075,2400)(4575,2400)(4575,2100)
\path(3675,2400)(3675,1800)
\path(3375,2400)(3375,1800)
\path(3975,2400)(3975,1800)
\path(4275,2400)(4275,1800)
\put(4800,0){\makebox(0,0)[lb]{{\SetFigFont{8}{9.6}{\rmdefault}{\mddefault}{\updefault}$k$}}}
\path(3075,300)(4875,300)
\path(4755.000,270.000)(4875.000,300.000)(4755.000,330.000)
\put(2850,2400){\makebox(0,0)[lb]{{\SetFigFont{8}{9.6}{\rmdefault}{\mddefault}{\updefault}$l$}}}
\path(5775,1800)(7275,1800)(7275,300)
\path(5775,1500)(7275,1500)
\path(5775,1200)(7275,1200)
\path(5775,900)(7275,900)
\path(5775,600)(7275,600)
\path(6075,1800)(6075,300)
\path(6375,1800)(6375,300)
\path(6675,1800)(6675,300)
\path(6975,1800)(6975,300)
\path(5775,300)(5775,900)
\path(5775,900)(5775,2700)
\path(5805.000,2580.000)(5775.000,2700.000)(5745.000,2580.000)
\path(5775,2100)(7275,2100)(7275,1800)
\path(5775,2400)(7275,2400)(7275,2100)
\path(6375,2400)(6375,1800)
\path(6075,2400)(6075,1800)
\path(6675,2400)(6675,1800)
\path(6975,2400)(6975,1800)
\put(7500,0){\makebox(0,0)[lb]{{\SetFigFont{8}{9.6}{\rmdefault}{\mddefault}{\updefault}$k$}}}
\path(5775,300)(7575,300)
\path(7455.000,270.000)(7575.000,300.000)(7455.000,330.000)
\put(5550,2400){\makebox(0,0)[lb]{{\SetFigFont{8}{9.6}{\rmdefault}{\mddefault}{\updefault}$l$}}}
\path(8475,1800)(9975,1800)(9975,300)
\path(8475,1500)(9975,1500)
\path(8475,1200)(9975,1200)
\path(8475,900)(9975,900)
\path(8475,600)(9975,600)
\path(8775,1800)(8775,300)
\path(9075,1800)(9075,300)
\path(9375,1800)(9375,300)
\path(9675,1800)(9675,300)
\path(8475,300)(8475,900)
\path(8475,900)(8475,2700)
\path(8505.000,2580.000)(8475.000,2700.000)(8445.000,2580.000)
\path(8475,2100)(9975,2100)(9975,1800)
\path(8475,2400)(9975,2400)(9975,2100)
\path(9075,2400)(9075,1800)
\path(8775,2400)(8775,1800)
\path(9375,2400)(9375,1800)
\path(9675,2400)(9675,1800)
\put(10200,0){\makebox(0,0)[lb]{{\SetFigFont{8}{9.6}{\rmdefault}{\mddefault}{\updefault}$k$}}}
\path(8475,300)(10275,300)
\path(10155.000,270.000)(10275.000,300.000)(10155.000,330.000)
\put(8250,2400){\makebox(0,0)[lb]{{\SetFigFont{8}{9.6}{\rmdefault}{\mddefault}{\updefault}$l$}}}
\path(11175,1800)(12675,1800)(12675,300)
\path(11175,1500)(12675,1500)
\path(11175,1200)(12675,1200)
\path(11175,900)(12675,900)
\path(11175,600)(12675,600)
\path(11475,1800)(11475,300)
\path(11775,1800)(11775,300)
\path(12075,1800)(12075,300)
\path(12375,1800)(12375,300)
\path(11175,300)(11175,900)
\path(11175,900)(11175,2700)
\path(11205.000,2580.000)(11175.000,2700.000)(11145.000,2580.000)
\path(11175,2100)(12675,2100)(12675,1800)
\path(11175,2400)(12675,2400)(12675,2100)
\path(11775,2400)(11775,1800)
\path(11475,2400)(11475,1800)
\path(12075,2400)(12075,1800)
\path(12375,2400)(12375,1800)
\put(12900,0){\makebox(0,0)[lb]{{\SetFigFont{8}{9.6}{\rmdefault}{\mddefault}{\updefault}$k$}}}
\path(11175,300)(12975,300)
\path(12855.000,270.000)(12975.000,300.000)(12855.000,330.000)
\put(10950,2400){\makebox(0,0)[lb]{{\SetFigFont{8}{9.6}{\rmdefault}{\mddefault}{\updefault}$l$}}}
\put(5880,1245){\makebox(0,0)[lb]{{\SetFigFont{8}{9.6}{\rmdefault}{\mddefault}{\updefault}1}}}
\put(6165,1545){\makebox(0,0)[lb]{{\SetFigFont{8}{9.6}{\rmdefault}{\mddefault}{\updefault}1}}}
\put(6450,1245){\makebox(0,0)[lb]{{\SetFigFont{8}{9.6}{\rmdefault}{\mddefault}{\updefault}1}}}
\put(8565,1230){\makebox(0,0)[lb]{{\SetFigFont{8}{9.6}{\rmdefault}{\mddefault}{\updefault}1}}}
\put(8565,915){\makebox(0,0)[lb]{{\SetFigFont{8}{9.6}{\rmdefault}{\mddefault}{\updefault}1}}}
\put(8865,930){\makebox(0,0)[lb]{{\SetFigFont{8}{9.6}{\rmdefault}{\mddefault}{\updefault}1}}}
\put(9450,915){\makebox(0,0)[lb]{{\SetFigFont{8}{9.6}{\rmdefault}{\mddefault}{\updefault}1}}}
\put(11265,930){\makebox(0,0)[lb]{{\SetFigFont{8}{9.6}{\rmdefault}{\mddefault}{\updefault}1}}}
\put(11580,645){\makebox(0,0)[lb]{{\SetFigFont{8}{9.6}{\rmdefault}{\mddefault}{\updefault}1}}}
\put(12465,675){\makebox(0,0)[lb]{{\SetFigFont{8}{9.6}{\rmdefault}{\mddefault}{\updefault}1}}}
\put(11850,645){\makebox(0,0)[lb]{{\SetFigFont{8}{9.6}{\rmdefault}{\mddefault}{\updefault}1}}}
\put(11565,930){\makebox(0,0)[lb]{{\SetFigFont{8}{9.6}{\rmdefault}{\mddefault}{\updefault}1}}}
\put(12165,930){\makebox(0,0)[lb]{{\SetFigFont{8}{9.6}{\rmdefault}{\mddefault}{\updefault}1}}}
\put(9180,1230){\makebox(0,0)[lb]{{\SetFigFont{8}{9.6}{\rmdefault}{\mddefault}{\updefault}1}}}
\put(480,45){\makebox(0,0)[lb]{{\SetFigFont{6}{7.2}{\rmdefault}{\mddefault}{\updefault}0}}}
\put(750,45){\makebox(0,0)[lb]{{\SetFigFont{5}{6.0}{\rmdefault}{\mddefault}{\updefault}1}}}
\put(1065,45){\makebox(0,0)[lb]{{\SetFigFont{5}{6.0}{\rmdefault}{\mddefault}{\updefault}2}}}
\put(1350,45){\makebox(0,0)[lb]{{\SetFigFont{5}{6.0}{\rmdefault}{\mddefault}{\updefault}3}}}
\put(1665,30){\makebox(0,0)[lb]{{\SetFigFont{5}{6.0}{\rmdefault}{\mddefault}{\updefault}4}}}
\put(465,1845){\makebox(0,0)[lb]{{\SetFigFont{8}{9.6}{\rmdefault}{\mddefault}{\updefault}1}}}
\put(3150,1860){\makebox(0,0)[lb]{{\SetFigFont{8}{9.6}{\rmdefault}{\mddefault}{\updefault}1}}}
\put(3495,1545){\makebox(0,0)[lb]{{\SetFigFont{8}{9.6}{\rmdefault}{\mddefault}{\updefault}1}}}
\put(120,1875){\makebox(0,0)[lb]{{\SetFigFont{5}{6.0}{\rmdefault}{\mddefault}{\updefault}0}}}
\put(30,1605){\makebox(0,0)[lb]{{\SetFigFont{5}{6.0}{\rmdefault}{\mddefault}{\updefault}-1}}}
\put(0,1275){\makebox(0,0)[lb]{{\SetFigFont{5}{6.0}{\rmdefault}{\mddefault}{\updefault}-2}}}
\put(15,960){\makebox(0,0)[lb]{{\SetFigFont{5}{6.0}{\rmdefault}{\mddefault}{\updefault}-3}}}
\put(15,675){\makebox(0,0)[lb]{{\SetFigFont{5}{6.0}{\rmdefault}{\mddefault}{\updefault}-4}}}
\path(375,1800)(1875,1800)(1875,300)
\end{picture}
}
\end{center}
\caption{$E_1$-pages of the spectral sequence for $j=-1,-3,-5,-7,-9$}
\label{fig:trefE1}
\end{figure}
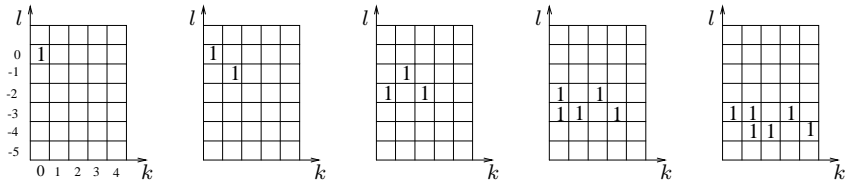

Using the computations in Example \ref{exe:trefoil} this gives 
$E_2$-pages as shown in Figure \ref{fig:trefE2}.
\begin{figure}[h]
\begin{center}
\setlength{\unitlength}{0.00033333in}
\begingroup\makeatletter\ifx\SetFigFont\undefined%
\gdef\SetFigFont#1#2#3#4#5{%
  \reset@font\fontsize{#1}{#2pt}%
  \fontfamily{#3}\fontseries{#4}\fontshape{#5}%
  \selectfont}%
\fi\endgroup%
{\renewcommand{\dashlinestretch}{30}
\begin{picture}(13053,2727)(0,-10)
\put(0,360){\makebox(0,0)[lb]{{\SetFigFont{5}{6.0}{\rmdefault}{\mddefault}{\updefault}-5}}}
\path(375,1500)(1875,1500)
\path(375,1200)(1875,1200)
\path(375,900)(1875,900)
\path(375,600)(1875,600)
\path(675,1800)(675,300)
\path(975,1800)(975,300)
\path(1275,1800)(1275,300)
\path(1575,1800)(1575,300)
\path(375,300)(375,900)
\path(375,900)(375,2700)
\path(405.000,2580.000)(375.000,2700.000)(345.000,2580.000)
\path(375,2100)(1875,2100)(1875,1800)
\path(375,2400)(1875,2400)(1875,2100)
\path(975,2400)(975,1800)
\path(675,2400)(675,1800)
\path(1275,2400)(1275,1800)
\path(1575,2400)(1575,1800)
\put(2100,0){\makebox(0,0)[lb]{{\SetFigFont{8}{9.6}{\rmdefault}{\mddefault}{\updefault}$k$}}}
\path(375,300)(2175,300)
\path(2055.000,270.000)(2175.000,300.000)(2055.000,330.000)
\put(150,2400){\makebox(0,0)[lb]{{\SetFigFont{8}{9.6}{\rmdefault}{\mddefault}{\updefault}$l$}}}
\path(3075,1800)(4575,1800)(4575,300)
\path(3075,1500)(4575,1500)
\path(3075,1200)(4575,1200)
\path(3075,900)(4575,900)
\path(3075,600)(4575,600)
\path(3375,1800)(3375,300)
\path(3675,1800)(3675,300)
\path(3975,1800)(3975,300)
\path(4275,1800)(4275,300)
\path(3075,300)(3075,900)
\path(3075,900)(3075,2700)
\path(3105.000,2580.000)(3075.000,2700.000)(3045.000,2580.000)
\path(3075,2100)(4575,2100)(4575,1800)
\path(3075,2400)(4575,2400)(4575,2100)
\path(3675,2400)(3675,1800)
\path(3375,2400)(3375,1800)
\path(3975,2400)(3975,1800)
\path(4275,2400)(4275,1800)
\put(4800,0){\makebox(0,0)[lb]{{\SetFigFont{8}{9.6}{\rmdefault}{\mddefault}{\updefault}$k$}}}
\path(3075,300)(4875,300)
\path(4755.000,270.000)(4875.000,300.000)(4755.000,330.000)
\put(2850,2400){\makebox(0,0)[lb]{{\SetFigFont{8}{9.6}{\rmdefault}{\mddefault}{\updefault}$l$}}}
\path(5775,1800)(7275,1800)(7275,300)
\path(5775,1500)(7275,1500)
\path(5775,1200)(7275,1200)
\path(5775,900)(7275,900)
\path(5775,600)(7275,600)
\path(6075,1800)(6075,300)
\path(6375,1800)(6375,300)
\path(6675,1800)(6675,300)
\path(6975,1800)(6975,300)
\path(5775,300)(5775,900)
\path(5775,900)(5775,2700)
\path(5805.000,2580.000)(5775.000,2700.000)(5745.000,2580.000)
\path(5775,2100)(7275,2100)(7275,1800)
\path(5775,2400)(7275,2400)(7275,2100)
\path(6375,2400)(6375,1800)
\path(6075,2400)(6075,1800)
\path(6675,2400)(6675,1800)
\path(6975,2400)(6975,1800)
\put(7500,0){\makebox(0,0)[lb]{{\SetFigFont{8}{9.6}{\rmdefault}{\mddefault}{\updefault}$k$}}}
\path(5775,300)(7575,300)
\path(7455.000,270.000)(7575.000,300.000)(7455.000,330.000)
\put(5550,2400){\makebox(0,0)[lb]{{\SetFigFont{8}{9.6}{\rmdefault}{\mddefault}{\updefault}$l$}}}
\path(8475,1800)(9975,1800)(9975,300)
\path(8475,1500)(9975,1500)
\path(8475,1200)(9975,1200)
\path(8475,900)(9975,900)
\path(8475,600)(9975,600)
\path(8775,1800)(8775,300)
\path(9075,1800)(9075,300)
\path(9375,1800)(9375,300)
\path(9675,1800)(9675,300)
\path(8475,300)(8475,900)
\path(8475,900)(8475,2700)
\path(8505.000,2580.000)(8475.000,2700.000)(8445.000,2580.000)
\path(8475,2100)(9975,2100)(9975,1800)
\path(8475,2400)(9975,2400)(9975,2100)
\path(9075,2400)(9075,1800)
\path(8775,2400)(8775,1800)
\path(9375,2400)(9375,1800)
\path(9675,2400)(9675,1800)
\put(10200,0){\makebox(0,0)[lb]{{\SetFigFont{8}{9.6}{\rmdefault}{\mddefault}{\updefault}$k$}}}
\path(8475,300)(10275,300)
\path(10155.000,270.000)(10275.000,300.000)(10155.000,330.000)
\put(8250,2400){\makebox(0,0)[lb]{{\SetFigFont{8}{9.6}{\rmdefault}{\mddefault}{\updefault}$l$}}}
\path(11175,1800)(12675,1800)(12675,300)
\path(11175,1500)(12675,1500)
\path(11175,1200)(12675,1200)
\path(11175,900)(12675,900)
\path(11175,600)(12675,600)
\path(11475,1800)(11475,300)
\path(11775,1800)(11775,300)
\path(12075,1800)(12075,300)
\path(12375,1800)(12375,300)
\path(11175,300)(11175,900)
\path(11175,900)(11175,2700)
\path(11205.000,2580.000)(11175.000,2700.000)(11145.000,2580.000)
\path(11175,2100)(12675,2100)(12675,1800)
\path(11175,2400)(12675,2400)(12675,2100)
\path(11775,2400)(11775,1800)
\path(11475,2400)(11475,1800)
\path(12075,2400)(12075,1800)
\path(12375,2400)(12375,1800)
\put(12900,0){\makebox(0,0)[lb]{{\SetFigFont{8}{9.6}{\rmdefault}{\mddefault}{\updefault}$k$}}}
\path(11175,300)(12975,300)
\path(12855.000,270.000)(12975.000,300.000)(12855.000,330.000)
\put(10950,2400){\makebox(0,0)[lb]{{\SetFigFont{8}{9.6}{\rmdefault}{\mddefault}{\updefault}$l$}}}
\put(5880,1245){\makebox(0,0)[lb]{{\SetFigFont{8}{9.6}{\rmdefault}{\mddefault}{\updefault}1}}}
\put(6165,1545){\makebox(0,0)[lb]{{\SetFigFont{8}{9.6}{\rmdefault}{\mddefault}{\updefault}1}}}
\put(6450,1245){\makebox(0,0)[lb]{{\SetFigFont{8}{9.6}{\rmdefault}{\mddefault}{\updefault}1}}}
\put(8565,1230){\makebox(0,0)[lb]{{\SetFigFont{8}{9.6}{\rmdefault}{\mddefault}{\updefault}1}}}
\put(9450,915){\makebox(0,0)[lb]{{\SetFigFont{8}{9.6}{\rmdefault}{\mddefault}{\updefault}1}}}
\put(12465,675){\makebox(0,0)[lb]{{\SetFigFont{8}{9.6}{\rmdefault}{\mddefault}{\updefault}1}}}
\put(12165,930){\makebox(0,0)[lb]{{\SetFigFont{8}{9.6}{\rmdefault}{\mddefault}{\updefault}1}}}
\put(9180,1230){\makebox(0,0)[lb]{{\SetFigFont{8}{9.6}{\rmdefault}{\mddefault}{\updefault}1}}}
\put(465,60){\makebox(0,0)[lb]{{\SetFigFont{6}{7.2}{\rmdefault}{\mddefault}{\updefault}0}}}
\put(750,60){\makebox(0,0)[lb]{{\SetFigFont{5}{6.0}{\rmdefault}{\mddefault}{\updefault}1}}}
\put(1050,45){\makebox(0,0)[lb]{{\SetFigFont{5}{6.0}{\rmdefault}{\mddefault}{\updefault}2}}}
\put(1380,15){\makebox(0,0)[lb]{{\SetFigFont{5}{6.0}{\rmdefault}{\mddefault}{\updefault}3}}}
\put(1680,45){\makebox(0,0)[lb]{{\SetFigFont{5}{6.0}{\rmdefault}{\mddefault}{\updefault}4}}}
\put(450,1845){\makebox(0,0)[lb]{{\SetFigFont{8}{9.6}{\rmdefault}{\mddefault}{\updefault}1}}}
\put(3165,1845){\makebox(0,0)[lb]{{\SetFigFont{8}{9.6}{\rmdefault}{\mddefault}{\updefault}1}}}
\put(3480,1530){\makebox(0,0)[lb]{{\SetFigFont{8}{9.6}{\rmdefault}{\mddefault}{\updefault}1}}}
\put(105,1875){\makebox(0,0)[lb]{{\SetFigFont{5}{6.0}{\rmdefault}{\mddefault}{\updefault}0}}}
\put(60,1545){\makebox(0,0)[lb]{{\SetFigFont{5}{6.0}{\rmdefault}{\mddefault}{\updefault}-1}}}
\put(45,1245){\makebox(0,0)[lb]{{\SetFigFont{5}{6.0}{\rmdefault}{\mddefault}{\updefault}-2}}}
\put(30,960){\makebox(0,0)[lb]{{\SetFigFont{5}{6.0}{\rmdefault}{\mddefault}{\updefault}-3}}}
\put(15,660){\makebox(0,0)[lb]{{\SetFigFont{5}{6.0}{\rmdefault}{\mddefault}{\updefault}-4}}}
\path(375,1800)(1875,1800)(1875,300)
\end{picture}
}
\end{center}
\caption{$E_2$-pages of the spectral sequence for $j=-1,-3,-5,-7,-9$}
\label{fig:trefE2}
\end{figure}
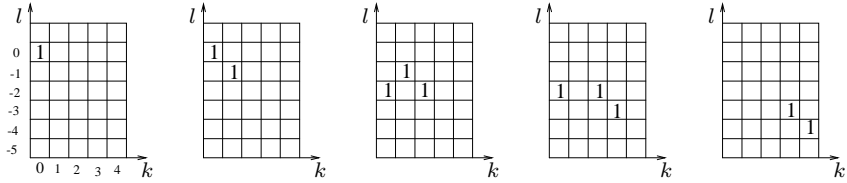

From this we can read off the Bar-Natan homology of the trefoil using
the fact that for dimensional reasons the spectral sequence collapses at $E_2$ and that $\bn i
j L \cong \bigoplus_{k+l = i} E_\infty^{k,l}$. Thus, for example, when
$j=-5$ we have $E_\infty$-term given in Figure \ref{fig:Einf}. 

\begin{figure}[h]
\begin{center}
\setlength{\unitlength}{0.00025000in}
\begingroup\makeatletter\ifx\SetFigFont\undefined%
\gdef\SetFigFont#1#2#3#4#5{%
  \reset@font\fontsize{#1}{#2pt}%
  \fontfamily{#3}\fontseries{#4}\fontshape{#5}%
  \selectfont}%
\fi\endgroup%
{\renewcommand{\dashlinestretch}{30}
\begin{picture}(5639,5697)(0,-10)
\put(3217,0){\makebox(0,0)[lb]{{\SetFigFont{7}{8.4}{\rmdefault}{\mddefault}{\updefault}$k+l=-2$}}}
\path(622,5070)(622,870)(3622,870)
\path(3622,870)(4222,870)
\blacken\path(4102.000,840.000)(4222.000,870.000)(4102.000,900.000)(4138.000,870.000)(4102.000,840.000)
\path(622,5070)(3622,5070)(3622,870)
\path(622,4470)(3622,4470)
\path(622,3870)(3622,3870)
\path(622,3270)(3622,3270)
\path(622,2670)(3622,2670)
\path(622,2070)(3622,2070)
\path(622,1470)(3622,1470)
\path(1222,5070)(1222,870)
\path(1822,5070)(1822,870)
\path(2422,5070)(2422,870)
\path(3022,5070)(3022,870)
\thicklines
\path(22,3870)(3622,270)
\path(22,5070)(4822,270)
\put(817,2805){\makebox(0,0)[lb]{{\SetFigFont{7}{8.4}{\rmdefault}{\mddefault}{\updefault}1}}}
\put(1417,3420){\makebox(0,0)[lb]{{\SetFigFont{7}{8.4}{\rmdefault}{\mddefault}{\updefault}1}}}
\put(2032,2820){\makebox(0,0)[lb]{{\SetFigFont{7}{8.4}{\rmdefault}{\mddefault}{\updefault}1}}}
\put(4747,630){\makebox(0,0)[lb]{{\SetFigFont{7}{8.4}{\rmdefault}{\mddefault}{\updefault}$k+l=0$}}}
\thinlines
\path(622,5070)(622,5670)
\blacken\path(652.000,5550.000)(622.000,5670.000)(592.000,5550.000)(622.000,5586.000)(652.000,5550.000)
\end{picture}
}
\end{center}
\caption{$E_\infty$-term for $j=-5$}
\label{fig:Einf}
\end{figure}
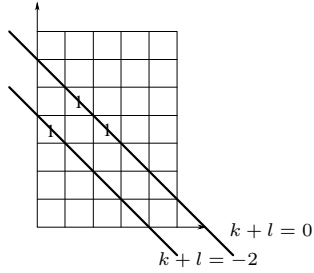

As there are only two non-trivial groups which lie on the lines $k+l=-2$ and $k+l=0$ we have  $\bn i j L
= 0$ unless $i=-2$ or $i=0$ in which case $\bn {-2} {-5} L = \FF$ and
$\bn 0 {-5} L = \FF \oplus \FF$.

Note that for $j<-9$ the
$E_\infty$-page will simply be a shifted version of the picture
above. Thus 
we can summarize the  graded Bar-Natan homology of the trefoil
in the table in Figure \ref{fig:BNtref}.
\begin{figure}[h]
\begin{center}
\setlength{\unitlength}{0.00033333in}
\begingroup\makeatletter\ifx\SetFigFont\undefined%
\gdef\SetFigFont#1#2#3#4#5{%
  \reset@font\fontsize{#1}{#2pt}%
  \fontfamily{#3}\fontseries{#4}\fontshape{#5}%
  \selectfont}%
\fi\endgroup%
{\renewcommand{\dashlinestretch}{30}
\begin{picture}(3387,4062)(0,-10)
\path(975,3012)(3375,3012)(3375,12)
	(975,12)(975,3012)
\path(975,2412)(3375,2412)
\path(975,1812)(3375,1812)
\path(3375,1212)(975,1212)
\path(975,612)(3375,612)
\path(1575,3012)(1575,12)
\path(2175,3012)(2175,12)
\path(2775,3012)(2775,12)
\put(2925,3162){\makebox(0,0)[lb]{{\SetFigFont{8}{9.6}{\rmdefault}{\mddefault}{\updefault}0}}}
\put(2400,3162){\makebox(0,0)[lb]{{\SetFigFont{8}{9.6}{\rmdefault}{\mddefault}{\updefault}-1}}}
\put(1200,3162){\makebox(0,0)[lb]{{\SetFigFont{8}{9.6}{\rmdefault}{\mddefault}{\updefault}-3}}}
\put(600,2562){\makebox(0,0)[lb]{{\SetFigFont{8}{9.6}{\rmdefault}{\mddefault}{\updefault}-1}}}
\put(600,1887){\makebox(0,0)[lb]{{\SetFigFont{8}{9.6}{\rmdefault}{\mddefault}{\updefault}-3}}}
\put(600,1362){\makebox(0,0)[lb]{{\SetFigFont{8}{9.6}{\rmdefault}{\mddefault}{\updefault}-5}}}
\put(600,762){\makebox(0,0)[lb]{{\SetFigFont{8}{9.6}{\rmdefault}{\mddefault}{\updefault}-7}}}
\put(1725,3162){\makebox(0,0)[lb]{{\SetFigFont{8}{9.6}{\rmdefault}{\mddefault}{\updefault}-2}}}
\put(1800,837){\makebox(0,0)[lb]{{\SetFigFont{8}{9.6}{\rmdefault}{\mddefault}{\updefault}1}}}
\put(1800,1437){\makebox(0,0)[lb]{{\SetFigFont{8}{9.6}{\rmdefault}{\mddefault}{\updefault}1}}}
\put(0,1437){\makebox(0,0)[lb]{{\SetFigFont{8}{9.6}{\rmdefault}{\mddefault}{\updefault}$j$}}}
\put(1950,3837){\makebox(0,0)[lb]{{\SetFigFont{8}{9.6}{\rmdefault}{\mddefault}{\updefault}$i$}}}
\put(50,237){\makebox(0,0)[lb]{{\SetFigFont{8}{9.6}{\rmdefault}{\mddefault}{\updefault}$j\leq9$}}}
\put(3000,2037){\makebox(0,0)[lb]{{\SetFigFont{8}{9.6}{\rmdefault}{\mddefault}{\updefault}2}}}
\put(3000,1437){\makebox(0,0)[lb]{{\SetFigFont{8}{9.6}{\rmdefault}{\mddefault}{\updefault}2}}}
\put(3000,237){\makebox(0,0)[lb]{{\SetFigFont{8}{9.6}{\rmdefault}{\mddefault}{\updefault}2}}}
\put(3000,837){\makebox(0,0)[lb]{{\SetFigFont{8}{9.6}{\rmdefault}{\mddefault}{\updefault}2}}}
\put(3000,2637){\makebox(0,0)[lb]{{\SetFigFont{8}{9.6}{\rmdefault}{\mddefault}{\updefault}1}}}
\end{picture}
}
\end{center}
\caption{Bar-Natan homology of the trefoil}
\label{fig:BNtref}
\end{figure}}
\end{exe}

\subsection{The stable range}\label{sec:stable}
It is possible to identify a {\em stable range} of Bar-Natan theory,
stable with respect to the $q$-grading. The stable range can be seen
in Example \ref{exe:tref} as the infinite tower of 2's. To define the stable range we use the following lemma.

\begin{lemma}
{\it There  exists $j_s\in\bZ$ such that the multiplication by $u$ map $u\colon \cch * j \ra \cch * {j-2}$ is an isomorphism for   $j\leq j_s$.}
\end{lemma}
\begin{proof}
The multiplication by $u$ map is clearly injective for all $j$. 

Let $j_s\in\bZ$ be such that $\cchb * j L =0$ for all
$j\leq j_s$. We claim that $u\colon \cch * j \ra \cch * {j-2}$ is
surjective for $j\leq j_s$. An element of $z \in\cch * {j-2}(L)$ can be
written in the form $z= \sum_\lambda u^{l_\lambda} y_\lambda$ where
$y_\lambda \in \cchb ** (L)$ and $l_\lambda\in \bN$. Since $q(z) = j-2$ we
have $j-2 = q( u^{l_\lambda} y_\lambda) = -2l_\lambda + q(y_\lambda)$ for
all $\lambda$ in the indexing set. Thus $l_\lambda -1 =
\frac{q(y_\lambda) - j}2 \geq 0$ since $q(y_\lambda) \geq j_s \geq
j$. Thus $\sum_\lambda u^{l_\lambda-1}y_\lambda$ is an element of $\cch *
j (L)$ which hits $z$ under multiplication by $u$.
\end{proof}

Since multiplication by $u$ is a chain map it follows that for $j\leq
j_s$ we have an isomorphism $u_*\colon\bn i {j} L \ra \bn i {j-2}
L$. Recalling that $\gamma$ is the number of components of $L$ modulo
2 we now define the singly graded {\em stable} Bar-Natan theory
as the direct limit
\[
\sbn i L = \lim (u_*\colon\bn i {2k+\gamma} L \ra \bn i {2k+\gamma -2}
L)
\]
Note that there is an isomorphism $\bn i j L \cong \sbn i L$ for each
$j\leq j_s$. 

As we shall see that stable theory is isomorphic to the filtered
($u=1$) Bar-Natan theory of Section~3. In fact the
spectral sequence in Subsection~4.1 for $j\leq j_s$ is
isomorphic to the spectral sequence in Subsection~3.2
for the filtered theory.

Recalling that an element of $\cch * {j}(L)$ can be written in the
form $\sum_\lambda u^{l_\lambda} y_\lambda$ we define $\eta\colon (\cch *
j (L),d) \ra (\scchb * (L),d)$ by $\eta(\sum_\lambda u^{l_\lambda} y_\lambda) =
y_\lambda$. This is clearly a chain map.

\begin{prop}
{\it Let $2j +\gamma \leq j_s$. Up to an overall shift in filtration degree
by $j$ the map $\eta\colon (\cch * {2j+\gamma} (L),d) \ra (\scchb *
(L),d)$ is an isomorphism of filtered complexes.}
\end{prop}
\begin{proof}
It is clear that $\eta$ is an injective chain map. To see it is
surjective let $\sum y_\lambda\in\scchb * (L)$. Then for each
$\lambda$ we have $q(y_\lambda)\geq j_s\geq 2j+\gamma$. Setting
$l_\lambda = \frac{q(y_\lambda) -2j -\gamma}2\geq 0$ we see
$\sum u^{l_\lambda} y_\lambda \in \cch * {2j+\gamma}$ and this
hits $\sum y_\lambda\in\scchb * (L)$ under $\eta$.

We now claim that $\eta$ maps $F^k\cch * {2j+\gamma} (L)$
isomorphically to $F^{j+k} \scchb * (L)$. Let $\sum
u^{l_\lambda} y_\lambda\in F^k\cch * {2j+\gamma} (L)$. Then by the
definition of the filtration $l_\lambda \geq k$ for all $\lambda$ in the
indexing set. We have $2j+\gamma = q(\sum u^{l_\lambda} y_\lambda)
= -2l_\lambda + q(y_\lambda)$ and so $q(y_\lambda) = 2j+\gamma +
2l_\lambda \geq 2(j+k) +\gamma$ showing that $\sum y_\lambda \in F^{j+k}
\scchb * (L)$.

Conversely, if $\sum y_\lambda \in F^{j+k} \scchb * (L)$ then
$q(y_\lambda)\geq 2(j+k) + \gamma$ and since $2j+\gamma = -2l_\lambda +
q(y_\lambda)$ we have $l_\lambda = \frac{q(y_\lambda) -2j -\gamma}2 \geq
\frac{2(j+k) + \gamma -2j - \gamma}2 = k$ and so $\sum
u^{l_\lambda} y_\lambda\in F^k\cch * {2j+\gamma} (L)$.
\end{proof}

It follows that $\eta$ induces a map of spectral sequences and up to an
shift of bi-degree $(j,-j)$ the $E_0$-pages of these two spectral sequences are
isomorphic, namely they are both given by the mod 2 Khovanov
chain complex. It follows that the higher pages of the two spectral
sequences are also isomorphic. Since the shift does not effect
the total degree we have the following corollary.

\begin{corollary}
\[
\sbn  i L \cong \ubn i 
\]
\end{corollary}

Thus the total dimension dim$\ubn *  = 2^k$ where $k$ is the number
of components and the homological degree of the generators is given by
the explicit calculation of Theorem \ref{thm:dim}.


\section{Applications}\label{sec:applications}

In this section we discuss two applications of the spectral sequences
defined above. We examine the appearance of ``shifted pawn moves'' in
Bar-Natan theory and describe the form of the mod 2 Khovanov
polynomial for thin knots.

\subsection{Tetris pieces and shifted pawn moves}

Bar-Natan has observed empirically that in many cases a tetris pieces in
$\FF$-Khovanov homology is replaced in his theory by a shifted pawn move. 

\begin{figure}[h]
\begin{center}
\setlength{\unitlength}{0.00033333in}
\begingroup\makeatletter\ifx\SetFigFont\undefined%
\gdef\SetFigFont#1#2#3#4#5{%
  \reset@font\fontsize{#1}{#2pt}%
  \fontfamily{#3}\fontseries{#4}\fontshape{#5}%
  \selectfont}%
\fi\endgroup%
{\renewcommand{\dashlinestretch}{30}
\begin{picture}(4973,2512)(0,-10)
\path(150,2485)(1350,2485)(1350,685)
	(150,685)(150,2485)
\path(750,2485)(750,685)
\path(150,1885)(1350,1885)
\path(150,1285)(1350,1285)
\put(900,1435){\makebox(0,0)[lb]{{\SetFigFont{8}{9.6}{\rmdefault}{\mddefault}{\updefault}1}}}
\put(900,2035){\makebox(0,0)[lb]{{\SetFigFont{8}{9.6}{\rmdefault}{\mddefault}{\updefault}1}}}
\path(3525,2485)(4725,2485)(4725,685)
	(3525,685)(3525,2485)
\path(4125,2485)(4125,685)
\path(3525,1885)(4725,1885)
\path(3525,1285)(4725,1285)
\put(4275,1435){\makebox(0,0)[lb]{{\SetFigFont{8}{9.6}{\rmdefault}{\mddefault}{\updefault}1}}}
\put(4275,2035){\makebox(0,0)[lb]{{\SetFigFont{8}{9.6}{\rmdefault}{\mddefault}{\updefault}1}}}
\put(2850,85){\makebox(0,0)[lb]{{\SetFigFont{8}{9.6}{\rmdefault}{\mddefault}{\updefault}Shifted
pawn move}}}
\put(300,1435){\makebox(0,0)[lb]{{\SetFigFont{8}{9.6}{\rmdefault}{\mddefault}{\updefault}1}}}
\put(300,835){\makebox(0,0)[lb]{{\SetFigFont{8}{9.6}{\rmdefault}{\mddefault}{\updefault}1}}}
\put(0,85){\makebox(0,0)[lb]{{\SetFigFont{8}{9.6}{\rmdefault}{\mddefault}{\updefault}Tetris piece}}}
\end{picture}
}
\end{center}
\end{figure}

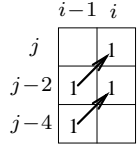
\begin{figure}[h,b,t]
\begin{center}
\setlength{\unitlength}{0.00033333in}
\begingroup\makeatletter\ifx\SetFigFont\undefined%
\gdef\SetFigFont#1#2#3#4#5{%
  \reset@font\fontsize{#1}{#2pt}%
  \fontfamily{#3}\fontseries{#4}\fontshape{#5}%
  \selectfont}%
\fi\endgroup%
{\renewcommand{\dashlinestretch}{30}
\begin{picture}(1662,2247)(0,-10)
\path(450,1812)(1650,1812)(1650,12)
	(450,12)(450,1812)
\path(1050,1812)(1050,12)
\path(450,1212)(1650,1212)
\path(450,612)(1650,612)
\thicklines
\path(750,912)(1200,1362)
\path(1072.721,1149.868)(1200.000,1362.000)(987.868,1234.721)
\path(750,312)(1200,762)
\path(1072.721,549.868)(1200.000,762.000)(987.868,634.721)
\put(1200,762){\makebox(0,0)[lb]{{\SetFigFont{8}{9.6}{\rmdefault}{\mddefault}{\updefault}1}}}
\put(600,762){\makebox(0,0)[lb]{{\SetFigFont{8}{9.6}{\rmdefault}{\mddefault}{\updefault}1}}}
\put(1200,1362){\makebox(0,0)[lb]{{\SetFigFont{8}{9.6}{\rmdefault}{\mddefault}{\updefault}1}}}
\put(600,162){\makebox(0,0)[lb]{{\SetFigFont{8}{9.6}{\rmdefault}{\mddefault}{\updefault}1}}}
\put(1250,1962){\makebox(0,0)[lb]{{\SetFigFont{8}{9.6}{\rmdefault}{\mddefault}{\updefault}$i$}}}
\put(450,1962){\makebox(0,0)[lb]{{\SetFigFont{8}{9.6}{\rmdefault}{\mddefault}{\updefault}$i\!-\!1$}}}
\put(0,1362){\makebox(0,0)[lb]{{\SetFigFont{8}{9.6}{\rmdefault}{\mddefault}{\updefault}$j$}}}
\put(-300,762){\makebox(0,0)[lb]{{\SetFigFont{8}{9.6}{\rmdefault}{\mddefault}{\updefault}$j\!-\!2$}}}
\put(-300,162){\makebox(0,0)[lb]{{\SetFigFont{8}{9.6}{\rmdefault}{\mddefault}{\updefault}$j\!-\!4$}}}
\end{picture}
}
\end{center}
\caption{Tetris piece with $\beta_*$ an isomorphism}
\label{fig:grtetris}
\end{figure}

Using the spectral sequence in Section~4.1 we can see
why this phenomenon occurs in certain circumstances. Suppose we have
the tetris piece in $\FF$-Khovanov homology with gradings as indicated
in Figure \ref{fig:grtetris}.

Suppose further that $\beta_*$ indicated by arrows are
isomorphisms. Then we
claim that this tetris piece contributes a shifted pawn move.

\begin{figure}[h]
\begin{center}
\setlength{\unitlength}{0.00041667in}
\begingroup\makeatletter\ifx\SetFigFont\undefined%
\gdef\SetFigFont#1#2#3#4#5{%
  \reset@font\fontsize{#1}{#2pt}%
  \fontfamily{#3}\fontseries{#4}\fontshape{#5}%
  \selectfont}%
\fi\endgroup%
{\renewcommand{\dashlinestretch}{30}
\begin{picture}(8434,2235)(0,-10)
\path(540,315)(2340,315)
\path(2220.000,285.000)(2340.000,315.000)(2220.000,345.000)
\path(540,315)(540,2115)
\path(570.000,1995.000)(540.000,2115.000)(510.000,1995.000)
\path(540,1815)(2040,1815)(2040,315)
\path(540,1515)(2040,1515)
\path(540,1215)(2040,1215)
\path(540,915)(2040,915)
\path(540,615)(2040,615)
\path(840,1815)(840,315)
\path(1140,1815)(1140,315)
\path(1440,1815)(1440,315)
\path(1740,1815)(1740,315)
\path(3570,315)(5370,315)
\path(5250.000,285.000)(5370.000,315.000)(5250.000,345.000)
\path(3570,315)(3570,2115)
\path(3600.000,1995.000)(3570.000,2115.000)(3540.000,1995.000)
\path(3570,1815)(5070,1815)(5070,315)
\path(3570,1515)(5070,1515)
\path(3570,1215)(5070,1215)
\path(3570,915)(5070,915)
\path(3570,615)(5070,615)
\path(3870,1815)(3870,315)
\path(4170,1815)(4170,315)
\path(4470,1815)(4470,315)
\path(4770,1815)(4770,315)
\path(6570,345)(8370,345)
\path(8250.000,315.000)(8370.000,345.000)(8250.000,375.000)
\path(6570,345)(6570,2145)
\path(6600.000,2025.000)(6570.000,2145.000)(6540.000,2025.000)
\path(6570,1845)(8070,1845)(8070,345)
\path(6570,1545)(8070,1545)
\path(6570,1245)(8070,1245)
\path(6570,945)(8070,945)
\path(6570,645)(8070,645)
\path(6870,1845)(6870,345)
\path(7170,1845)(7170,345)
\path(7470,1845)(7470,345)
\path(7770,1845)(7770,345)
\put(390,1965){\makebox(0,0)[lb]{{\SetFigFont{11}{13.2}{\rmdefault}{\mddefault}{\updefault}$l$}}}
\put(2265,15){\makebox(0,0)[lb]{{\SetFigFont{11}{13.2}{\rmdefault}{\mddefault}{\updefault}$k$}}}
\put(615,15){\makebox(0,0)[lb]{{\SetFigFont{8}{13.2}{\rmdefault}{\mddefault}{\updefault}$0$}}}
\put(250,1275){\makebox(0,0)[lb]{{\SetFigFont{8}{13.2}{\rmdefault}{\mddefault}{\updefault}$i$}}}
\put(600,1245){\makebox(0,0)[lb]{{\SetFigFont{11}{13.2}{\rmdefault}{\mddefault}{\updefault}1}}}
\put(3390,1965){\makebox(0,0)[lb]{{\SetFigFont{11}{13.2}{\rmdefault}{\mddefault}{\updefault}$l$}}}
\put(6390,1980){\makebox(0,0)[lb]{{\SetFigFont{11}{13.2}{\rmdefault}{\mddefault}{\updefault}$l$}}}
\put(5295,0){\makebox(0,0)[lb]{{\SetFigFont{11}{13.2}{\rmdefault}{\mddefault}{\updefault}$k$}}}
\put(8265,45){\makebox(0,0)[lb]{{\SetFigFont{11}{13.2}{\rmdefault}{\mddefault}{\updefault}$k$}}}
\put(3310,1275){\makebox(0,0)[lb]{{\SetFigFont{8}{13.2}{\rmdefault}{\mddefault}{\updefault}$i$}}}
\put(2900,960){\makebox(0,0)[lb]{{\SetFigFont{8}{13.2}{\rmdefault}{\mddefault}{\updefault}$i\!-\!1$}}}
\put(5915,960){\makebox(0,0)[lb]{{\SetFigFont{8}{13.2}{\rmdefault}{\mddefault}{\updefault}$i\!-\!1$}}}
\put(5915,660){\makebox(0,0)[lb]{{\SetFigFont{8}{13.2}{\rmdefault}{\mddefault}{\updefault}$i\!-\!2$}}}
\put(3630,1245){\makebox(0,0)[lb]{{\SetFigFont{11}{13.2}{\rmdefault}{\mddefault}{\updefault}1}}}
\put(3645,945){\makebox(0,0)[lb]{{\SetFigFont{11}{13.2}{\rmdefault}{\mddefault}{\updefault}1}}}
\put(3945,945){\makebox(0,0)[lb]{{\SetFigFont{11}{13.2}{\rmdefault}{\mddefault}{\updefault}1}}}
\put(6630,975){\makebox(0,0)[lb]{{\SetFigFont{11}{13.2}{\rmdefault}{\mddefault}{\updefault}1}}}
\put(6945,990){\makebox(0,0)[lb]{{\SetFigFont{11}{13.2}{\rmdefault}{\mddefault}{\updefault}1}}}
\put(6975,675){\makebox(0,0)[lb]{{\SetFigFont{11}{13.2}{\rmdefault}{\mddefault}{\updefault}1}}}
\put(7260,675){\makebox(0,0)[lb]{{\SetFigFont{11}{13.2}{\rmdefault}{\mddefault}{\updefault}1}}}
\end{picture}
}
\end{center}
\caption{$E_1$-pages of the spectral sequence for $j$, $j-2$ and $j-4$}
\label{fig:tetE1}
\end{figure}
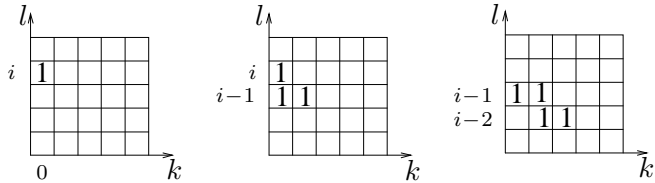

The $E_1$-pages for the spectral sequences
corresponding to $j, j-2$ and $j-4$ are given in Figure
\ref{fig:tetE1}.

\begin{figure}[h]
\begin{center}
\setlength{\unitlength}{0.00041667in}
\begingroup\makeatletter\ifx\SetFigFont\undefined%
\gdef\SetFigFont#1#2#3#4#5{%
  \reset@font\fontsize{#1}{#2pt}%
  \fontfamily{#3}\fontseries{#4}\fontshape{#5}%
  \selectfont}%
\fi\endgroup%
{\renewcommand{\dashlinestretch}{30}
\begin{picture}(8434,2235)(0,-10)
\path(540,315)(2340,315)
\path(2220.000,285.000)(2340.000,315.000)(2220.000,345.000)
\path(540,315)(540,2115)
\path(570.000,1995.000)(540.000,2115.000)(510.000,1995.000)
\path(540,1815)(2040,1815)(2040,315)
\path(540,1515)(2040,1515)
\path(540,1215)(2040,1215)
\path(540,915)(2040,915)
\path(540,615)(2040,615)
\path(840,1815)(840,315)
\path(1140,1815)(1140,315)
\path(1440,1815)(1440,315)
\path(1740,1815)(1740,315)
\path(3570,315)(5370,315)
\path(5250.000,285.000)(5370.000,315.000)(5250.000,345.000)
\path(3570,315)(3570,2115)
\path(3600.000,1995.000)(3570.000,2115.000)(3540.000,1995.000)
\path(3570,1815)(5070,1815)(5070,315)
\path(3570,1515)(5070,1515)
\path(3570,1215)(5070,1215)
\path(3570,915)(5070,915)
\path(3570,615)(5070,615)
\path(3870,1815)(3870,315)
\path(4170,1815)(4170,315)
\path(4470,1815)(4470,315)
\path(4770,1815)(4770,315)
\path(6570,345)(8370,345)
\path(8250.000,315.000)(8370.000,345.000)(8250.000,375.000)
\path(6570,345)(6570,2145)
\path(6600.000,2025.000)(6570.000,2145.000)(6540.000,2025.000)
\path(6570,1845)(8070,1845)(8070,345)
\path(6570,1545)(8070,1545)
\path(6570,1245)(8070,1245)
\path(6570,945)(8070,945)
\path(6570,645)(8070,645)
\path(6870,1845)(6870,345)
\path(7170,1845)(7170,345)
\path(7470,1845)(7470,345)
\path(7770,1845)(7770,345)
\put(390,1965){\makebox(0,0)[lb]{{\SetFigFont{11}{13.2}{\rmdefault}{\mddefault}{\updefault}$l$}}}
\put(2265,15){\makebox(0,0)[lb]{{\SetFigFont{11}{13.2}{\rmdefault}{\mddefault}{\updefault}$k$}}}
\put(615,15){\makebox(0,0)[lb]{{\SetFigFont{8}{13.2}{\rmdefault}{\mddefault}{\updefault}$0$}}}
\put(250,1275){\makebox(0,0)[lb]{{\SetFigFont{8}{13.2}{\rmdefault}{\mddefault}{\updefault}$i$}}}
\put(600,1245){\makebox(0,0)[lb]{{\SetFigFont{11}{13.2}{\rmdefault}{\mddefault}{\updefault}1}}}
\put(3390,1965){\makebox(0,0)[lb]{{\SetFigFont{11}{13.2}{\rmdefault}{\mddefault}{\updefault}$l$}}}
\put(6390,1980){\makebox(0,0)[lb]{{\SetFigFont{11}{13.2}{\rmdefault}{\mddefault}{\updefault}$l$}}}
\put(5295,0){\makebox(0,0)[lb]{{\SetFigFont{11}{13.2}{\rmdefault}{\mddefault}{\updefault}$k$}}}
\put(8265,45){\makebox(0,0)[lb]{{\SetFigFont{11}{13.2}{\rmdefault}{\mddefault}{\updefault}$k$}}}
\put(3310,1275){\makebox(0,0)[lb]{{\SetFigFont{8}{13.2}{\rmdefault}{\mddefault}{\updefault}$i$}}}
\put(2900,960){\makebox(0,0)[lb]{{\SetFigFont{8}{13.2}{\rmdefault}{\mddefault}{\updefault}$i\!-\!1$}}}
\put(5915,960){\makebox(0,0)[lb]{{\SetFigFont{8}{13.2}{\rmdefault}{\mddefault}{\updefault}$i\!-\!1$}}}
\put(5915,660){\makebox(0,0)[lb]{{\SetFigFont{8}{13.2}{\rmdefault}{\mddefault}{\updefault}$i\!-\!2$}}}
\put(3630,1245){\makebox(0,0)[lb]{{\SetFigFont{11}{13.2}{\rmdefault}{\mddefault}{\updefault}1}}}
\end{picture}
}
\end{center}
\caption{$E_2$-pages of the spectral sequence for $j$, $j-2$ and $j-4$}
\label{fig:tetE2}
\end{figure}
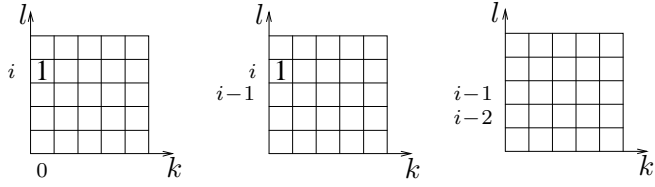

Using the fact that $d_1=\beta_*$  is an
isomorphism on the relevant groups we get $E_2$-pages as in
Figure \ref{fig:tetE2}.

From the $E_2$ page for $j$ we see we have a contribution of $\FF$ to
$\bn {i} {j} L$.  From the $E_2$ page for $j-2$ we get an $\FF$ in
$\bn {i} {j-2} L = \FF$ and no contribution to $\bn {i-1} {j-2} L$ and
finally from the $E_2$ page for $j-4$ we see there is no contribution
at all. Assuming further that these groups we have identified at $E_2$
survive to $E_\infty$ we get the situation summarized in Figure \ref{fig:grpawn}.

\begin{figure}[h]
\begin{center}
\setlength{\unitlength}{0.00033333in}
\begingroup\makeatletter\ifx\SetFigFont\undefined%
\gdef\SetFigFont#1#2#3#4#5{%
  \reset@font\fontsize{#1}{#2pt}%
  \fontfamily{#3}\fontseries{#4}\fontshape{#5}%
  \selectfont}%
\fi\endgroup%
{\renewcommand{\dashlinestretch}{30}
\begin{picture}(1662,2247)(0,-10)
\path(450,1812)(1650,1812)(1650,12)
	(450,12)(450,1812)
\path(1050,1812)(1050,12)
\path(450,1212)(1650,1212)
\path(450,612)(1650,612)
\put(1200,762){\makebox(0,0)[lb]{{\SetFigFont{8}{9.6}{\rmdefault}{\mddefault}{\updefault}1}}}
\put(1200,1362){\makebox(0,0)[lb]{{\SetFigFont{8}{9.6}{\rmdefault}{\mddefault}{\updefault}1}}}
\put(1250,1962){\makebox(0,0)[lb]{{\SetFigFont{8}{9.6}{\rmdefault}{\mddefault}{\updefault}$i$}}}
\put(450,1962){\makebox(0,0)[lb]{{\SetFigFont{8}{9.6}{\rmdefault}{\mddefault}{\updefault}$i\!-\!1$}}}
\put(0,1362){\makebox(0,0)[lb]{{\SetFigFont{8}{9.6}{\rmdefault}{\mddefault}{\updefault}$j$}}}
\put(-300,762){\makebox(0,0)[lb]{{\SetFigFont{8}{9.6}{\rmdefault}{\mddefault}{\updefault}$j\!-\!2$}}}
\put(-300,162){\makebox(0,0)[lb]{{\SetFigFont{8}{9.6}{\rmdefault}{\mddefault}{\updefault}$j\!-\!4$}}}
\end{picture}
}
\end{center}
\caption{}
\label{fig:grpawn}
\end{figure}

\subsection{The mod 2 Khovanov polynomial for thin knots.}
In this case, by reverse engineering, knowledge of the explicit
calculation of the filtered theory and the spectral sequence allows
us to deduce the form of the mod 2 Khovanov polynomial for thin knots.

\begin{prop}\label{prop:es}
{\it 
Let $L$ be an oriented link. If the spectral sequence of Section~3.2
collapses at $E_2$ then for each $k$ there is a sequence of groups
{\small
\[
\xymatrix{ \ar[r] & \kh {i-1}{2k+\gamma + 2(i-1)} L
  \ar[r]^-{\beta_*} & \kh {i}{2k+\gamma + 2i} L
  \ar[r]^-{\beta_*} & \kh {i+1}{2k+\gamma + 2(i+1)} L
  \ar[r]^-{\beta_*} &
}
\]
}
which is exact at $\kh {i}{2k+\gamma + 2i} L $ whenever $\ubn i  = 0$.
}
\end{prop}

\begin{proof}
If $\ubn i = 0$ then $\kk i {2j_s + \gamma +2p} L = 0$ for all
$p$. This follows from the explicit calculation of the filtered theory
and by the assumption that the spectral sequence collapses at
$E_2$. In other words if $\ubn i = 0$ then $\kk i j L = 0$ for all $j$, from which it
follows that
\[
\frac{\Ker (\beta_*\colon \kh i j L \ra \kh {i+1} {j+2} L)}
{\Image (\beta_*\colon \kh {i-1} {j-2} L \ra \kh {i} {j} L)}= 0.
\]
\end{proof}

For a knot $L$ the above sequence is exact except possibly at the 0'th group.
There is a similar sequence for the spectral sequence of reduced theory.

\vspace*{6pt}

\begin{thm}\label{thm:thin}
{\it If $L$ is an $\FF$H-thin knot with homology concentrated on
diagonals $j=s-1+2i$ and $j=s+1+2i$ then there exists a polynomial
$Kh^\prime (L)$ such that 
\[
Kh_{\FF}(L) = q^{s-1}(1+q^2)(1+(1+tq^2)Kh^\prime (L))
\]
where  $Kh^\prime (L)$ is a polynomial in $tq^2$.}
\end{thm}

\vspace*{6pt}

\begin{proof}
We work with reduced theory which has homology concentrated on the single diagonal $j=s+2i$. For dimensional reasons the spectral sequence for reduced theory collapses at the $E_2$-page and by the reduced version of Proposition \ref{prop:es} we have a sequence
\[
{\small
\xymatrix{
 \ar[r] & \rkh {i-1}{s+2(i-1)} L \ar[r]^-{\beta_*^{i-1}} & \rkh
  {i}{s+2i} L \ar[r]^-{\beta_*^{i}} & \rkh {i+1}{s+2(i+1)} L
    \ar[r]^-{\beta_*^{i+1}} &  
}
}
\]
which is exact except when $i=0$. For $i\neq 0$ we have
\begin{eqnarray*}
\dim \rkh i {s+2i} L & = &\dim \Image \beta_*^i + \dim \Ker
\beta_*^i\\
 & = &  \dim \Ker \beta_*^{i+1} + \dim \Ker \beta_*^i\\
 & =: & K_{i+1} + K_{i}.
\end{eqnarray*}

When $i=0$ the deviation from exactness is given by knowledge of the
(reduced) filtered theory, which tells us there is one additional
generator.  Thus the reduced mod 2 Khovanov polynomial
$\widetilde{Kh}_{\FF}(L)$ is of the form
\begin{eqnarray*}
q^{s} + \sum_i t^iq^{s+2i} (K_{i+1} + K_{i}) & = &
q^{s} + \sum_i q^{s} (t^{i-1}q^{2(i-1)} + t^iq^{2i})K_{i}\\
 & = &  q^{s}(1+  \sum_i  (t^{i-1}q^{2(i-1)} + t^iq^{2i})K_{i})\\
 & = &  q^{s}(1+  (1+tq^2) \sum_i  t^{i-1}q^{2(i-1)}K_{i})\\
& = & q^{s}(1+  (1+tq^2) \sum_i  (tq^2)^{i-1}K_{i}).
\end{eqnarray*}

Using results of Shumakovitch \cite{shumakovitch} (see also
\cite{rasmussen2}) in the mod 2 case we have $\kh ** L \cong \rkh ** L \otimes A$ from which it follows that
\[
Kh_{\FF} (L) = \widetilde{Kh}_{\FF}(L)(q+q^{-1})=
(q+q^{-1})q^{s}(1+(1+tq^2) \sum_i (tq^2)^{i-1} K_i)
\]
as claimed.
\end{proof}

The formula in the proposition above first appeared in \cite{barnatan1} where  
$Kh^\prime$ is the polynomial calculated in that paper for rational
Khovanov theory\footnote{It appeared
in the context that in an earlier version of his paper he conjectured
that for any prime knot there existed an integer $s$ making the above
hold true. However, knot $8_{19}$ provided a counter example.}. Here we do not claim that the polynomial $Kh^\prime (L)$ of the
proposition is the same as Bar-Natan's.

\section{Proof of Proposition \ref{prop:commute}}\label{app}

As we have already mentioned the proof of Proposition
\ref{prop:commute}  is  very similar to the proof given by Lee in
\cite{lee}. For convenience we include the details here. This
section, however, may be seen as an exposition of her work in a
slightly modified context. The reader worried about any sign
discrepancies with Lee's work should remember we are working over
$\FF$.

We will use the following notation and conventions given a link
diagram $D$. Firstly, throughout this subsection we are going to omit all
shifts of the various complexes involved which will unclutter the
notation a little. If $\epsilon$ is a string of $k$ 0's and 1's then we can
consider the subcube (of the cube of smoothings) consisting of all
smoothings where the last $k$ entries are given by $\epsilon$. Denote
the associated complex by
$\cC(D(*\epsilon))$. Given $k$ then $\cC
(D)$ may be decomposed as (a graded group but not as a complex)
\[
\cC(D) \cong \bigoplus_\epsilon \cC(D(*\epsilon))
\]
where $\epsilon$ runs over all strings of 0's and 1's of length
$k$. With respect to this decomposition the differential $\partial$
(or the map $\beta$) may be decomposed as well. Supposing $v\in
\cC(D(*\epsilon))$ then we write
\[
\partial (v) = \dd \epsilon {} (v) + \sum_{\epsilon^\prime} \dd
\epsilon {\epsilon^\prime} (v)
\]
where $\epsilon^\prime$ has one more 1 than $\epsilon$ and $\dd
\epsilon {\epsilon^\prime} \colon \cd \epsilon \ra \cd
	 {\epsilon^\prime}$ and $\dd \epsilon {} \colon \cd \epsilon
	 \ra \cd \epsilon$ are induced by the differential $\partial$.

\subsubsection{Reidemeister I positive twist}
We refer to Figure \ref{fig:RI} below.

\begin{figure}[h]
\centerline{
\psfig{figure=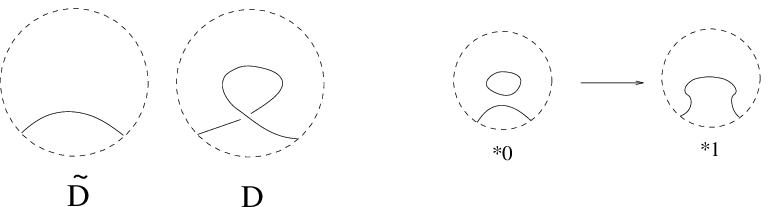}
}
\caption{Reidemeister I}
\label{fig:RI}
\end{figure}

There are isomorphisms of complexes
\begin{equation}\label{R1isos}
\cC (\tilde D ) \cong \cd 1 \;\;\;\; \mbox{and} \;\;\;\; \cd 1 \ot A
\cong \cd 0
\end{equation}
and these will be used throughout.

Khovanov's quasi-isomorphism for this move is defined as follows.
$\cC(D)$ decomposes into a direct sum of subcomplexes
\[
\cC(D) \cong X_1 \oplus X_2
\]
where $X_2$ is contractible. These two subcomplexes are defined by
\begin{eqnarray*}
X_1 & = & \Ker (\dd 0 1) \subset \cd 0\\
X_2 & = & \{ a\ot 1 + b \mid a,b\in \cd 1 \} \subset \cd 0 \oplus \cd 1
\end{eqnarray*}
Note that using the isomorphism
(\ref{R1isos}) any element of $\cd 0$ (and hence of $X_1$) can be
written in the form $a\ot 1 + b\ot x$ for $a,b\in \cd 1$. There is a
map $\rho \colon X_1 \ra \cd 1$ defined by 
\[
\rho ( a\ot 1 + b\ot x) = b
\]
and the composite
\[
\xymatrix{ 
\phi \colon \cC (D) \cong X_1 \oplus X_2   \ar[r]^-{pr}  & X_1
  \ar[r]^-{\rho} & \cd 1 \cong \cC (\tilde D)
}
\]
is a quasi-isomorphism (an isomorphism on homology).

Now let $ a\ot 1 + b\ot x $ be a cycle in $X_1$. In a moment we will
need the following identity:
\begin{equation}\label{R1ident}
\dd 0 1 (\bb 1 {} a \ot 1) + \dd 0 1 (\bb 1 {} b \ot 1 ) + \dd 1 {}
\bb 0 1 ( b\ot x) = 0.
\end{equation}
This can be shown to hold by recalling that $\beta$ is a map of
complexes so $\partial \beta =
\beta\partial$ and then expanding the right hand side of $0 = \beta (0) =
\beta (\partial (a\ot 1 + b\ot x)) = \partial(\beta(a\ot 1 + b\ot
x))$.

Now we compute
\begin{align*}
\beta (a\ot 1 + & b\ot x)  =  \bb 0 {} (a\ot 1) + \bb 0 {} (b\ot x) +
\bb 0 1 (a\ot 1) + \bb 0 1 (b\ot x)\\
 & =  \bb 1 {}a \ot 1 + \bb 1 {}b \ot x + \bb 01 (b\ot x)\\
& =  (\bb 1 {} a + \dd 1{}\bb 01 (b\ot x)) \ot 1 + \bb 1 {} b \ot x +
\dd 1{}\bb 01 (b\ot x)\ot 1 + \bb 01 (b\ot x).
\end{align*}
Using (\ref{R1ident}) on can show that $\dd 01 ((\bb 1 {} a + \dd
1{}\bb 01 (b\ot x)) \ot 1 + \bb 1 {} b \ot x )= 0$ and hence $(\bb 1 {}
a + \dd 1{}\bb 01 (b\ot x)) \ot 1 + \bb 1 {} b \ot x \in X_1$. Note
also that $\dd 1{}\bb 01 (b\ot x)\ot 1 + \bb 01 (b\ot x) \in
X_2$. Thus we have
\[
\phi\beta (a\ot 1 + b\ot x) = \rho ((\bb 1 {} a + \dd
1{}\bb 01 (b\ot x)) \ot 1 + \bb 1 {} b \ot x) = \bb 1{} b.
\]
More easily we see
\[
\beta\phi (a\ot 1 + b\ot x) = \beta (b) = \bb 1 {} b
\]
and so $\beta$ and $\phi$ commute.


\subsubsection{Reidemeister II}
We refer to Figure \ref{fig:RII} below.

\begin{figure}[h]
\centerline{
\psfig{figure=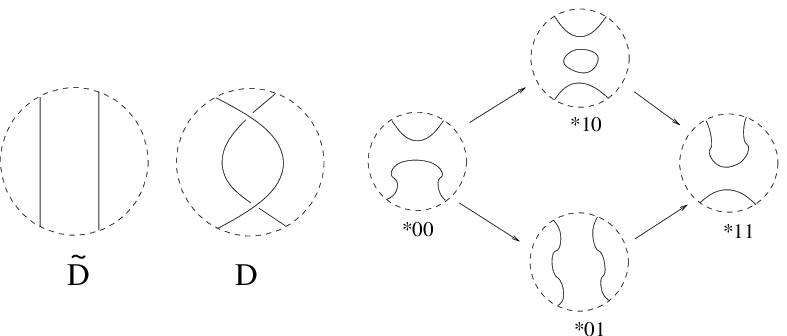}
}
\caption{Reidemeister II}
\label{fig:RII}
\end{figure}

There are isomorphisms of complexes
\begin{equation}\label{R2isos}
\cC (\tilde D ) \cong \cd {01} \;\;\;\; \mbox{and} \;\;\;\; \cd
    {11}\ot A\cong \cd {10}.
\end{equation}
There is a decomposition
\[\cC (D) = X_1 \oplus X_2 \oplus X_3\]
where $X_2$ and $X_3$ are contractible. The complex $X_1$ is defined by
\begin{eqnarray*}
X_1 & = & \{a + \dd {01}{11} a \ot 1\mid a \in \cd {01}\}  \subset \cd
{01} \oplus \cd {10}
\end{eqnarray*}
There is a map $\rho\colon \cd
{01} \ra X_1$ defined by
\[
\rho (a) = a + \dd {01}{11} a \ot 1
\]
where the right hand side uses the second isomorphism in
\ref{R2isos}. The composite
\[
\xymatrix{ 
\cC(\tilde D) \cong \cd {01}   \ar[r]^-{\rho}  & X_1
  \ar@{^{(}->}[r] & \cC (D)
}
\]
is a quasi-isomorphism. 

Let $a$ be a cycle in $\cd {01}$ so $\dd {01}{}a = 0$. Regarding $a\in
\cC(D)$ we use the equality $\partial \beta + \beta\partial = 0$ to
deduce
\begin{equation}\label{R2ident1}
\bb{11}{} \dd{01}{11} a + \dd{01}{11}\bb{01}{} a = \dd{11}{} \bb{01}{11}a.
\end{equation}
Also note that for any $z\in \cd {11}$ we have
\begin{eqnarray}
\label{R2ident2} \dd {10}{11} (z\ot 1) & = & z,\\
\label{R2ident3} \bb {10}{11} (z\ot 1) & = & 0,\\
\label{R2ident4} \bb {10}{} (z\ot 1) & = & \bb {11}{} z \ot 1.
\end{eqnarray}

Now 
 \begin{align*}
\beta \rho (a) = \beta (a + \dd {01}{11} a \ot 1) & = \bb{01}{} a + \bb
      {01}{11} a + \bb{10}{} (\dd {01}{11} a \ot 1) + \bb{10}{11} (\dd
      {01}{11} a \ot 1)\\
& = \bb{01}{} a + \bb{01}{11} a + \bb{11}{}\dd{01}{11} a \ot 1
      \;\;\;\;\mbox{by (\ref{R2ident3}) and (\ref{R2ident4})}
 \end{align*}
Hence we have
\begin{align*}
\beta \rho (a) + \rho\beta (a) & = (\bb{01}{} a + \bb{01}{11} a +
\bb{11}{}\dd{01}{11} a \ot 1) + (\bb{01}{} a + \dd{01}{11}\bb{01}{} a
\ot 1) \\
 & =  \bb{01}{11} a +
\bb{11}{}\dd{01}{11} a \ot 1 +  \dd{01}{11}\bb{01}{} a
\ot 1\\
 & =  \bb{01}{11} a + \dd{11}{} \bb{01}{11} a\otimes 1 \;\;\;\;\;\;\mbox{ by
  (\ref{R2ident1})}\\
& = \dd{10}{11}(\bb{01}{11} a \ot 1) + \dd{01}{} (\bb{01}{11} a \ot 1)
\;\;\;\;\;\;\mbox{ by
  (\ref{R2ident2})}\\
 & = \partial (\bb{01}{11} a \ot 1).
\end{align*}
Thus $\beta$ and $\rho$ commute up to boundaries.

\subsubsection{Reidemeister III}
We refer to Figure \ref{fig:RIII} below.

\begin{figure}[h,t,b]
\centerline{
\psfig{figure=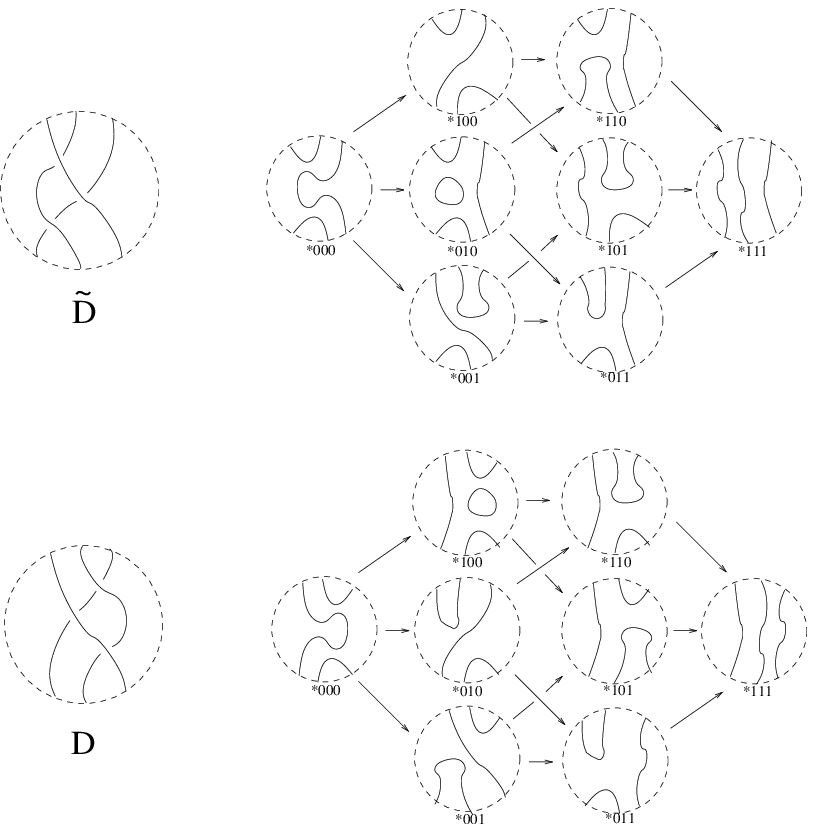}
}
\caption{Reidemeister III}
\label{fig:RIII}
\end{figure}

There are isomorphisms of complexes
\begin{equation}\label{R3isos1}
\cdt {110} \otimes A \cong \cdt{010} \;\;\;\; \mbox{and} \;\;\;\;\;
\cd{110} \ot A \cong \cd {100}
\end{equation}
and also
\begin{equation}\label{R3isos2}
\xymatrix{
\cdt {100} \ar[r]_f^{\cong} & \cd{010} & \;\;\mbox{and} \;\; & 
\cdt{1}  \ar[r]_g^{\cong} & \cd {1}
}
\end{equation}

There are decompositions
\[\cC (\tilde D) = \tilde X_1 \oplus \tilde X_2 \oplus \tilde
X_3\;\;\;\; \mbox{and} \;\;\;\;\;\cC (D) =  X_1 \oplus  X_2 \oplus  X_3\]
where $X_2$, $X_3$, $\tilde X_2$ and $\tilde X_3$ are
contractible. The complexes we are required to know in detail are
\begin{eqnarray*}
\tilde X_1 & = & \{a + \dd {100}{110} a \ot 1 + b \mid a \in \cdt
       {100}, b\in \cdt 1\}  
\\
\tilde X_3 & = & \{ a\ot 1 + \partial (b\ot 1) \mid a, b \in \cdt
       {110}\}\\
X_1  & = & \{a + \dd {010}{110} a \ot 1 + b \mid a \in \cd
       {110}, b\in \cd 1\}  
\end{eqnarray*}
There is a map $\rho\colon \tilde X_1 \ra X_1$ defined by
\[
\rho (a + \dd {100}{110} a \ot 1 + b) = f(a) + \dd {010}{110} f(a) \ot 1 + g(b) 
\]
and the composite
\[
\xymatrix{ 
\phi\colon \cC(\tilde D) \ar[r]^-{pr} & \tilde X_1 \ar[r]^-{\rho}  & X_1
  \ar@{^{(}->}[r] & \cC (D)
}
\]
is a quasi-isomorphism. 

Let $a + \dd {100}{110} a \ot 1 + b $ be a cycle in $\tilde X_1$. In
particular this implies $\dd{100}{}a = 0$. Looking at the $\cdt {110}$
component of $(\partial \beta + \beta\partial)(a) = 0$ 
we get the following identity.
\begin{equation}\label{R3ident1}
\bb{110}{}\dd{100}{110} a  = \dd{100}{110}\bb{100}{} a +
\dd{110}{}\bb{100}{110}a
\end{equation}

Also for $z\in\cdt {110}$ we have
\begin{equation}\label{R3ident3}
\bb{010}{} (z\ot 1) = \bb{110}{} z \ot 1
\end{equation}

In order to calculate $\phi\beta$ we need to first get an expression
for $\beta (a + \dd {100}{110} a \ot 1 + b) $ in terms of the
decomposition of $\cC (\tilde D)$ above.
\begin{align*}
&\beta(a +   \dd {100}{110} a \ot  1 + b) \\
& =  \bb{100}{}a \!+\! \bb{100}{101}a \!+\! \bb{100}{110}a \!+\!
    \bb{010}{}(\dd{100}{110}a\ot 1) \!+\! \bb{010}{110} (\dd{100}{110}
    a\ot 1)\!+\! \bb{010}{011}(\dd{100}{110}a\ot 1) \!+ \!\bb{1}{} b\\
& =  \bb{100}{}a + \bb{100}{101}a + \bb{100}{110}a +
    \bb{010}{}(\dd{100}{110}a\ot 1)  + \bb{1}{} b\\
& =  \bb{100}{}a + \bb{100}{101}a + \bb{100}{110}a +
    \bb{110}{}\dd{100}{110}a\ot 1  + \bb{1}{} b \;\;\;\;\mbox{by
      (\ref{R3ident3})}\\
& =  \bb{100}{}a + \bb{100}{101}a + \bb{100}{110}a +
   \dd{100}{110}\bb{100}{} a \otimes 1+ \dd{110}{}\bb{100}{110}a\otimes 1  + \bb{1}{} b 
\;\;\;\;\mbox{by (\ref{R3ident1})}\\
 & = \bb{100}{}a + \dd{100}{110}\bb{100}{} a \ot 1+ \bb{100}{101}a + 
\bb{1}{} b + \dd{010}{011}(\bb{100}{110}a\ot 1) + \partial (\bb{100}{110}a \ot 1)   
\end{align*}
The first five summands give an element of $\tilde X_1$ and the last
an element of $\tilde X_3$. 

A similar computation gives that 
\begin{align*}
\beta (f(a) + & \dd {010}{110} f(a) \ot 1 + g(b))\\
&  = \bb{010}{}f(a) +
\bb{010}{011}f(a) + \dd{010}{110}\bb{010}{}f(a)\ot 1 +
\dd{100}{101}(\bb{010}{110}f(a) \ot 1) + \bb{1}{} g(b) \\
& \;\;\;\;\; + \partial
(\bb{010}{110}f(a)\ot 1).
\end{align*}

From the above we can now see that 
\begin{align*}
\phi\beta (a + & \dd {100}{110} a \ot 1 + b)) \\
& = f(\bb{100}{}a) + \dd{010}{110}f(\bb{100}{} a) \ot 1 + g(\bb{100}{101}a + 
\bb{1}{} b + \dd{010}{011}(\bb{100}{110}a\ot 1)).
\end{align*}
Keeping the definitions of $\beta$ and $\partial$ in mind, looking at
Figure \ref{fig:RIII} we see that 
\[
f(\bb{100}{}a) = \bb{010}{}f(a), \;\;\;\;\;\;\; f(\bb{100}{110}a) = \bb{010}{110} f(a), 
\;\;\;\;\;\; g(\bb 1{} b) = \bb 1 {} g(b),
\]
\[
g(\bb{100}{101}a) = \dd{100}{101}
(\bb{010}{110} f(a) \ot 1), \;\;\;\;\; g(\dd{010}{011}(\bb{100}{110} a
\ot 1 ))= \bb{010}{011}f(a),
\]
so that 
\begin{align*}
\phi\beta (a + & \dd {100}{110} a \ot 1 + b)) \\
& = 
\bb{010}{}f(a) +
\bb{010}{011}f(a) + \dd{010}{110}\bb{010}{}f(a)\ot 1 +
\dd{100}{101}(\bb{010}{110}f(a) \ot 1) + \bb{1}{} g(b).
\end{align*}

Hence we finally get
\begin{align*}
(\phi \beta + \beta \phi)(a +  \dd {100}{110} a & \ot 1 + b) \\
&  = 
\phi\beta(a + \dd {100}{110} a \ot 1 + b) + \beta (f(a) + \dd
	 {010}{110} f(a) \ot 1 + g(b))\\
& =    \partial
(\bb{010}{110}f(a)\ot 1)
\end{align*}

and so $\beta$ and $\phi$ commute up to boundaries.

\vspace*{1cm}

\noindent {\bf Acknowledgements}
The author was supported by the European Commission through a Marie
Curie fellowship and thanks the Institut de Recherche Math\'ematique
Avanc\'ee in Strasbourg for their hospitality. Thanks also to
J. Rasmussen for comments on an early version of this work and to the
referee for further suggestions for improvement.


\end{document}